\documentclass[12pt, a4paper]{article}

\usepackage[utf8]{inputenc}
\usepackage{titlesec}
\usepackage[english]{babel}
\usepackage{graphicx}
\usepackage{float}
\usepackage{fancyhdr}
\usepackage[matrix,arrow,curve]{xy}
\usepackage{pstricks} 
\usepackage{amsmath,amsfonts,verbatim,afterpage,theorem,euscript,mathrsfs,amssymb}
\usepackage{amsfonts}
\usepackage{amssymb}
\usepackage{array}
\usepackage{dsfont}
\usepackage{hyperref}
\usepackage{authblk}

\newtheorem{Definition}{Definition}[section]
\newtheorem{Proposition}{Proposition}[section]

\newtheorem{Lemma}{Lemma}[section]
\newtheorem{Theorem}{Theorem}

\def \ve{\textbf{e}}
\def \vu{\textbf{u}}
\def \vb{\textbf{b}}
\def \vv{\textbf{v}}
\def \vw{\textbf{w}}

\newcommand\vN{{\bf \nabla}}

\def \pv{{\bf{Proof.}}~}

\newcommand \Endproof{\hfill $\diamond$}

\begin{document}

\title{Weighted energy estimates for the incompressible Navier-Stokes equations and applications to axisymmetric solutions without swirl}
\author{Pedro Gabriel Fern\'andez-Dalgo\footnote{Universit\'e Paris-Saclay, CNRS, Univ Evry, Laboratoire de Math\'ematiques et Mod\'elisation d'Evry, 91037, Evry-Courcouronnes, France} \footnote{e-mail: pedro.fernandez@univ-evry.fr} ,  Pierre Gilles Lemari\'e-Rieusset\footnote{Universit\'e Paris-Saclay, CNRS, Univ Evry, Laboratoire de Math\'ematiques et Mod\'elisation d'Evry, 91037, Evry-Courcouronnes, 
 France} \footnote{e-mail: pierregilles.lemarierieusset@univ-evry.fr}}
\date{}\maketitle

\maketitle

\begin{abstract}
We consider a family of weights which permit to generalize the Leray procedure to obtain weak suitable solutions of the 3D incompressible Navier--Stokes equations with initial data in weighted $L^2$ spaces. Our principal result concerns the existence of regular global solutions when the initial velocity is an axisymmetric vector field without swirl  such that both the initial velocity and its vorticity belong to $L^2 ( (1+ r^2)^{-\frac{\gamma}{2}} dx ) $, with $r= \sqrt{x_1^2 + x_2^2}$ and $\gamma \in (0, 2) $.
 \end{abstract}
 

\section{Introduction}

In 1934, Leray \cite{Le34} proved global existence of weak solutions for the 3D incompressible Navier--Stokes equations     \begin{equation*}  (NS) 
        \left\{ \begin{matrix} 
            \partial_t \vu= \Delta \vu  -(\vu\cdot \vN)\vu- \vN p 
            \cr \cr \vN\cdot \vu=0, \phantom{space space} \vu(0,.)=\vu_0
        \end{matrix}\right.
    \end{equation*}  in the case of a fluid filling the whole space whose initial velocity $\vu_0$  is in $L^2$.  Leray's strategy is to regularize the initial value and to mollify the non-linearity through convolution with a bump function:  let $\theta_\epsilon (x)=\frac 1{\epsilon^3}\theta(\frac x \epsilon )$, where $\theta\in\mathcal{D}(\mathbb{R}^3)$, $\theta$ is non-negative and radially decreasing and $\int\theta\, dx=1$; the mollified equations are then
      \begin{equation*}  (NS_\epsilon) 
        \left\{ \begin{matrix} 
            \partial_t \vu_\epsilon  = \Delta \vu_\epsilon  -( (\theta_\epsilon*\vu_\epsilon) \cdot \vN) \vu_\epsilon - \vN p_\epsilon  
            \cr \cr \vN \cdot \vu_\epsilon =0, \phantom{space space} \vu_\epsilon(0,.)=\theta_\epsilon*\vu_{0}.
        \end{matrix}\right.
    \end{equation*}
Standard methods give existence of a smooth solution on an interval $[0,T_\epsilon]$ where $T_\epsilon\approx \epsilon^3 \| \theta_\epsilon*\vu_0\|_2^{-2}$. Then, the energy equality
$$ \|\vu_\epsilon(t,.)\|_2^2+2\int_0^t \|\vN\otimes\vu_\epsilon\|_2^2\, ds=\|\theta_\epsilon*\vu_0\|_2^2$$
allows one  to extend the existence time and to get a global solution $\vu_\epsilon$; moreover, the same energy equality allows one to use a compactness argument and to get a subsequence $\vu_{\epsilon_k}$ that converges to a solution $\vu$ of the Navier--Stokes equations (NS) which satisfies the energy \emph{inequality}
$$ \|\vu(t,.)\|_2^2+2\int_0^t \|\vN\otimes\vu\|_2^2\, ds\leq \|\vu_0\|_2^2.$$
Weak solutions of equations (NS) that satisfy this energy inequality are called Leray solutions.
 
 There are many ways to extend Leray's results to settings where $\vu_0$ has infinite energy. A natural one is based on a splitting $\vu=\vv+\vw$ where $\vv$ satisfies an equation  $$\partial_t\vv=\Delta \vv+F(\vv), \text{ where }\nabla\cdot  F(\vv)=0,$$that is easy to solve and $\vw$ satisfies   perturbed Navier--Stokes equations $$\partial_t\vw+\vw\cdot\nabla\vw= \Delta \vw -\nabla q-\vv\cdot\nabla\vw-\vw\cdot\nabla\vv -\vv\cdot\nabla\vv-F(\vv)$$ for which Leray's formalism still holds. For instance, if $\vu_0\in L^p$ with $2\leq p<3$, Calder\' on \cite{Cal90, LRX} splits $\vu_0$ into $\vv_0+\vw_0$ where $\vv_0$ is a divergence-free vector field which is small in $L^3$ and $\vw_0$ belongs to $L^2$. Then, $\vv$ is the mild solution of the Navier--Stokes problem with $\vv_0$ as initial value. Another recent example is the way Seregin and \v Sver\'ak \cite{SerS} deal with global weak solutions for large initial values in $L^3$, by splitting $\vu$ into $\vv+\vw$, where  $\vv=e^{t\Delta}\vu_0$ and $\vw$ which has a finite energy.
 
 Another way to extend Leray's method is to consider weighted energy inequalities in $L^2(\Phi\, dx)$, or similarly energy inequalities for   $\sqrt\Phi \vu$. $\sqrt \Phi \vu$ is solution of  \begin{equation*}   
        \left\{ \begin{matrix} 
            \partial_t (\sqrt\Phi \vu)= \sqrt\Phi \Delta \vu  -\sqrt\Phi(\vu\cdot \vN)\vu- \sqrt\Phi\vN p 
            \cr \cr \vN\cdot \vu=0, \phantom{space space} \sqrt\Phi \vu(0,.)=\sqrt\Phi\vu_0
        \end{matrix}\right.
    \end{equation*}
    The problem is that the non-linear part of the equation $ -\sqrt\Phi(\vu\cdot \vN)\vu- \sqrt\Phi\vN p 
$ will then contribute to the energy balance, in contrast to the case of the Leray method. More precisely, we may write
 $$ -\sqrt\Phi(\vu\cdot \vN)\vu  = -\vu\cdot \nabla(\sqrt\Phi \vu)- \sqrt\Phi \vu \cdot (\sqrt\Phi\vu\cdot \nabla\frac 1{\sqrt\Phi}).
$$ 
 The advection term $\vu\cdot \nabla(\sqrt\Phi \vu)$ corresponds to a transport  by a divergence-free vector field and will not contribute to the energy; in order to control the impact of $\sqrt\Phi \vu\cdot (\sqrt\Phi\vu\cdot \nabla\frac 1{\sqrt\Phi})$ on $\sqrt\Phi \vu$, it is natural to assume that $\nabla \frac 1{\sqrt\Phi}$ is bounded; assuming that $\Phi$ is positive, we find that $\frac 1{\sqrt{\Phi(\frac{\lambda x}{|x|})}}\leq \frac 1{\sqrt{\Phi(0)}}+ C \lambda$, and thus that $\frac 1{1+\vert x\vert^2}\leq C \Phi(x)$.

Recently, Bradshaw, Kukavica  and Tsai \cite{BKT}  and Fern\'andez-Dalgo and Lemari\'e-Rieusset \cite{PF_PG}  used Leray's procedure to find a global weak solution to the equations (NS) when $\vu_0$ is no longer assumed to have finite energy but only to satisfy the weaker assumption $$\int \vert \vu_0(x)\vert^2\ \frac{dx}{1+\vert x\vert^2}<+\infty.$$
The solutions then satisfy, for every finite positive $T$,
$$\sup_{0\leq t\leq T} \int \vert \vu(t,x)\vert^2\ \frac{dx}{1+\vert x\vert^2} +\int_0^T \int \vert \vN\otimes\vu(t,x)\vert^2\ \frac{dx}{1+\vert x\vert^2}  <+\infty.$$

For the proof, a precise description of the presssure is needed, as it interfers as well in the energy balance; This point has been discussed in \cite{BTpressure} and \cite{PF_PG2}. The scheme of proof of existence of such weak solutions can easily be generalized to equations which behave  like the Navier--Stokes equations, for instance the magneto-hydrodynamic equations  \cite{PF_OJ, PF_OJ2}. 

An application of solutions in $L^2(\frac 1{1+\vert x\vert^2}dx)$ is given in \cite{PF_PG}: a simple proof of existence of discretely self-similar solutions to the Navier--Stokes problem when the initial value is locally  square integrable (and discretely homogeneous: $\lambda \vu_0(\lambda x)=\vu_0(x)$ for some $\lambda>1$). This existence was first proved by Chae and Wolf \cite{CW18}, and Bradshaw and Tsai \cite{BT19}, as a generalization of the result of Jia and \v Sver\'ak \cite{JS14}  for a regular homogeneous initial value ($\lambda \vu_0(\lambda x)=\vu_0(x)$ for every $\lambda>1$), see  \cite{LR16} for the case of locally square-integrable homogeneous initial value. If $\vu_0$ is homogeneous and locally square-integrable, then it belongs to $L^2_{\rm uloc}$;  thus, the proof of Jia and \v Sver\'ak relied on the control of weak solutions in the space $L^2_{\rm uloc}$ of uniformly locally square integrable vector fields, following the theory developed in \cite{LR02}.  If $\vu_0$ is discretely homogeneous and locally square-integrable, then it may fail to belong  to $L^2_{\rm uloc}$, but it belongs to $L^2(\frac 1{1+\vert x\vert^\gamma}\, dx)$ for $\gamma>1$.\\

Whereas the cases of finite energy and of infinite energy sound very similar, this similarity breaks down when we consider higher regularity.

When we consider solutions in function spaces with decaying weights, the growth of solutions can be amplified by the non-linearities. The authors in \cite{BKT} and \cite{PF_PG} used the transport structure of the non-linearity $ ( \vu \cdot \vN ) \vu $ to get good controls for the velocity in some weighted spaces.
When dealing with derivatives of the velocity, one loses the transport structure of non-linearities. The problem comes from the stretching term $\omega\cdot\nabla\vu$ in the equations for the vorticity
 $$\partial_t \omega = \Delta \omega + ( \omega \cdot \vN) \vu - (\vu \cdot \vN) \omega. $$ In the case when $\vu_0$ belongs to the classical Sobolev space $H^1$, for which local existence of a unique mild solution is known, this stretching term may potentially lead to blow-up in finite time, since it has a non-linear impact on the growth of $\|\omega\|_2$. There are two cases when this impact can be controlled: the case of 2D fluids (as the stretching term is equal to $0$) and the case of axisymmetric vector fields with no swirl \cite{La68, LMNP}. In the case of weighted estimates, one cannot even get local control of the size of the vorticity in $L^2(\Phi\, dx)$ in general, but we shall show that we have global existence of a weak solution such that $\|\sqrt\Phi \omega(t,.)\|_2$ remains  bounded on every bounded interval of times, when we work in 2D or when we consider  axisymmetric vector fields with no swirl  and  weights that depend  only on the distance to the symmetry axis.

\section{Main results.}
   We shall first prove global existence in the weighted $L^2$ setting, in dimension $d$ with $2\leq d\leq 4$ when the weight $\Phi$ satisfies some basic assumptions that allow the use of Leray's projection operator   and of energy estimates:

\begin{Definition}
An adapted  weight function $\Phi$ on $\mathbb{R}^d$ ($2\leq d\leq 4$) is a continuous  Lipschitz function $\Phi$ such that:
 \begin{itemize}
    \item $(\textbf{H}1)$ $0 < \Phi \leq  1$.    
    
    \item $(\textbf{H}2)$ There exists $C_1>0$ such that $|\vN \Phi | \leq C_1  \Phi ^{\frac{3}{2}} $
    
    \item $(\textbf{H}3)$ There exists $r \in (1, 2 ] $ such that  $\Phi^r \in \mathcal{A}_{r}$ (where $\mathcal{A}_r$ is the Muckenhoupt class of weights). In the case $d=4$, we require $r<2$ as well.

    \item $(\textbf{H}4)$ There exists $C_2>0$ such that  $ \Phi(x)\leq \Phi(\frac{x}{\lambda}) \leq C_2 \lambda^2 \Phi (x)$, for all $\lambda \geq 1$. 
    
\end{itemize}
\end{Definition}

Examples of adapted weights can easily be given by radial slowly  decaying functions:
\begin{itemize}
\item[$\bullet$] $d=2$, $\Phi(x)=\frac 1{(1+\vert x\vert)^\gamma}$ where $0\leq\gamma<2$
\item[$\bullet$] $d=3$ or $d=4$, $\Phi(x)=\frac 1{(1+\vert x\vert)^\gamma}$ where $0\leq\gamma\leq 2$
\item[$\bullet$] $d=3$, $\Phi(x)=\frac 1{(1+  r)^\gamma}$ where  $r=\sqrt{x_1^2+x_2^2}$ and $0\leq\gamma<2$.
\end{itemize}

  $\ $

  The following result concerns the existence of weak solutions belonging to a weighted $L^2$ space, where the weight permits to consider initial data with a weak decay at infinity.

\begin{Theorem}
\label{Thm1}
Let $d \in \{2, 3, 4 \}$. Consider a weight $\Phi$ satisfying $(\textbf{H}1)-(\textbf{H}4)$.
Let $\vu_0$ be a divergence free vector field, such that $\vu_0 $ belongs to $L^2(\Phi \,  dx, \mathbb{R}^d)$.   
Then, there exists a global solution $\vu$ of the problem
    \begin{equation*}  (NS) 
        \left\{ \begin{matrix} 
            \partial_t \vu= \Delta \vu  -(\vu\cdot \vN)\vu- \vN p 
            \cr \cr \vN\cdot \vu=0, \phantom{space space} \vu(0,.)=\vu_0
        \end{matrix}\right.
    \end{equation*}
such that
 \begin{itemize}
    \item $\vu$  belongs to $ L^\infty ((0,T), L^2 (\Phi dx) ) $ and $\vN \otimes\vu$ belongs to $L^2((0,T),L^2 (\Phi dx) )$, for all $T>0$,
    
    \item $ p = \sum_{1 \leq i,j \leq d} R_i R_j (u_i u_j)$,
     
    \item the map $t\in [0,+\infty)\mapsto \vu(t,.)$ is weakly continuous from $[0,+\infty)$ to $L^2 (\Phi \, dx) $, and  strongly continuous at $t=0$,
        
    \item For $d \in \{ 2, 3 \}$, $\vu$ satisfies the local energy inequality: there exists a locally finite non-negative measure $\mu$ such that
    \begin{equation*}
        \partial_t(\frac { | \vu | ^2}2)= \Delta (\frac { | \vu | ^2}2) - | \vN\otimes \vu | ^2 - \vN\cdot\left( \frac{ | \vu | ^2}2\vu\right)-\vN\cdot(p\vu)  -\mu,
    \end{equation*}
   and we have $\mu=0$ when $d=2$.
 \end{itemize}
\end{Theorem}

We observe that we do not prove the local energy inequality for the solutions in dimension $4$. We refer the papers \cite{DoGu}, \cite{WaWu} and \cite{Wu} for more information on suitable solutions in dimension 4.

If we consider the problem of higher regularity, the case of dimension $d=2$ is easy, while, in the case $d=3$, one must restrict the study to the case of axisymmetric flows with no swirl (to circumvent the stretching effect in the evolution of the vorticity).
\begin{Theorem} [Case $d=2$.]
\label{Thm2}
Let $ \Phi $ be a weight satisfying $(\textbf{H}1)-(\textbf{H}4) $.
Let $\vu_0$ be a divergence free   vector field, such that $\vu_0, \vN \otimes\vu_0 $ belong to $ L^2 (\Phi dx) $. Then there exists a global solution $\vu$ of the problem
    \begin{equation*}  (NS) 
        \left\{ \begin{matrix} 
            \partial_t \vu= \Delta \vu  -(\vu\cdot \vN)\vu- \vN p 
            \cr \cr \vN\cdot \vu=0, \phantom{space space} \vu(0,.)=\vu_0
        \end{matrix}\right.
    \end{equation*}
such that
 \begin{itemize}
    \item  $\vu$ and $\vN\otimes \vu $ belong to $ L^\infty ((0,T), L^2 (\Phi\, dx) ) $ and $ \Delta \vu $ belongs to $L^2((0,T),L^2 (\Phi \, dx) )$, for all $T>0$,
     
    \item the maps $t\in [0,+\infty)\mapsto \vu(t,.)$ and $t\in [0,+\infty)\mapsto \vN \otimes\vu(t,.)$ are weakly continuous from $[0,+\infty)$ to $L^2 (\Phi dx) $, and are strongly continuous at $t=0$.
         
 \end{itemize}
\end{Theorem}
\begin{Theorem} [Case $d=3$.]
\label{Thm3}
Let $ \Phi $ be a weight satisfying $(\textbf{H}1)-(\textbf{H}4) $.
Let $\vu_0$ be a divergence free axisymmetric vector field without swirl, such that $\vu_0, \vN\otimes \vu_0 $ belong to $ L^2 (\Phi \, dx) $. Assume moreover that $\Phi $ depends only on $r=\sqrt{x_1^2+x_2^2}$. Then there exists a time $T>0$, and a local solution $\vu$ on $(0,T)$ of the problem
    \begin{equation*}  (NS) 
        \left\{ \begin{matrix} 
            \partial_t \vu= \Delta \vu  -(\vu\cdot \vN)\vu- \vN p 
            \cr \cr \vN\cdot \vu=0, \phantom{space space} \vu(0,.)=\vu_0
        \end{matrix}\right.
    \end{equation*}
such that
 \begin{itemize}
    \item $\vu $ is axisymmetric without swirl, $\vu$ and $\vN\otimes \vu $ belong to $ L^\infty ((0,T), L^2 (\Phi\,  dx) ) $ and $ \Delta \vu $ belongs to $L^2((0,T),L^2 (\Phi\,  dx) )$,
     
    \item the maps $t\mapsto \vu(t,.)$ and $t\mapsto \vN \vu(t,.)$ are weakly continuous from $[0,T)$ to $L^2 (\Phi\,  dx) $, and are strongly continuous at $t=0$.
         
 \end{itemize}
\end{Theorem}

An extra condition on the weight permits to obtain a global existence result. Moreover, if the vorticity is more integrable at time $t=0$,  it will remain so in positive times. The next theorem precise these conditions on the weight.

\begin{Theorem} [Case $d=3$.]
\label{Thm4}
Let $ \Phi $ be a weight satisfying $(\textbf{H}1)-(\textbf{H}4) $.  Assume moreover that $\Phi $ depends only on $r=\sqrt{x_1^2+x_2^2}$.  Let $\Psi$ be another continuous weight (that  depends only on $r$)   such that  $\Phi \leq \Psi\leq 1$, $\Psi\in\mathcal{A}_2$ and       there exists $C_1>0$ such that $$ \vert \vN \Psi  \vert \leq C_1  \sqrt\Phi  \Psi \text{ and }  \vert\Delta \Psi\vert\leq C_1 \Phi\Psi.$$

Let $\vu_0$ be a divergence free axisymmetric vector field without swirl, such that $\vu_0,$ belongs to $ L^2 (\Phi dx) $ and $\vN\otimes \vu_0 $ belongs to $ L^2 (\Psi dx) $.   Then there exists a global solution $\vu$ of the problem
    \begin{equation*}  (NS) 
        \left\{ \begin{matrix} 
            \partial_t \vu= \Delta \vu  -(\vu\cdot \vN)\vu- \vN p 
            \cr \cr \vN\cdot \vu=0, \phantom{space space} \vu(0,.)=\vu_0
        \end{matrix}\right.
    \end{equation*}
such that
 \begin{itemize}
    \item $\vu $ is axisymmetric without swirl, $\vu$    belongs to $ L^\infty ((0,T), L^2 (\Phi\,  dx) ) $,  $\vN\otimes \vu $ belong to $ L^\infty ((0,T), L^2 (\Psi\,  dx) ) $ and $ \Delta \vu $ belongs to $L^2((0,T),L^2 (\Psi\,  dx) )$, for all $T>0$,
     
    \item the maps $t\in [0,+\infty)\mapsto \vu(t,.)$ and $t\in [0,+\infty)\mapsto \vN \otimes\vu(t,.)$ are weakly continuous from $[0,+\infty)$ to $L^2 (\Phi\,  dx) $ and to $L^2(\Psi\, dx)$ respectively, and are strongly continuous at $t=0$.
         
 \end{itemize}
\end{Theorem}

Example: we can take  $\Phi(x)=\frac 1{(1+  r)^\gamma}$  and $\Psi(x)=\frac 1{(1+  r^2)^{\delta/2}}$  with $0\leq\delta\leq\gamma<2$. Of course, $\Phi\approx \frac 1{(1+  r^2)^{\gamma/2}}$. The case $\delta<\gamma$ means that, if $\omega$ has a better decay at initial time, it will keep this beter decay at all times.

\section{Some lemmas on weights.}

Let us first recall the definition  of Muckenhoupt weights: for $1<q<+\infty$, a positive weight $ \Phi $  belongs to 
 $\mathcal{A}_q(\mathbb{R}^d) $ if and only if  \begin{equation}\label{muck}
  \sup_{x\in\mathbb{R}^d, \rho>0}  \left(\frac 1{\vert B(x,\rho)\vert}\int_{B(x,\rho)} \Phi \, dx\right)^{\frac{1}{q}} \left( \frac 1 {\vert B(x,\rho)\vert}\int_{B(x,\rho)} \Phi^{-\frac{1}{q-1}} \, dx \right)^{1-\frac{1}{q}}  <+\infty.\end{equation}
  We refer to the Chapter 9 in \cite{Gr09}.

  Due to the H\"older inequality, we have  $\mathcal{A}_q(\mathbb{R}^d) \subset \mathcal{A}_r(\mathbb{R}^d) $ if $q\leq r$.

 One easily cheks that $w_\gamma=\frac 1{(1+\vert x\vert)^\gamma}$ belongs to $\mathcal{A}_q(\mathbb{R}^d)$ if and only if $$ -d(q-1)<\gamma<d.$$
Thus, $\Phi=w_\gamma$ is an adapted weight if and only if $0\leq\gamma\leq 2$ and $\gamma<d$.

One may of course replace in inequality (\ref{muck}) the balls $B(x,\rho)$ by the cubes $Q(x,\rho)=]x_1-\rho,x_1+\rho[\times\dots\times]x_d-\rho,x_d+\rho[$. Thus, we can see that,
if $\Phi(x) = \Psi(x_1, x_2)$ and $1< q < + \infty $, then $\Phi \in \mathcal{A}_q(\mathbb{R}^3) $ if and only if $\Psi \in \mathcal{A}_q(\mathbb{R}^2) $. In particular, $\Phi(x)=\frac 1{(1+  r)^\gamma}$ is an adapted weight on $\mathbb{R}^3$ if and only if $0\leq\gamma<2$.

\begin{Lemma}\label{lem1} Let $\Phi$ satisfy $(\textbf{H}1)$ and $(\textbf{H} 2)$ and let $1\leq r<+\infty$. Then:
\\ a) $\sqrt \Phi f\in H^1$ if and only if $f\in L^2(\Phi\, dx)$ and $\vN f\in L^2(\Phi\, dx)$; moreover we have
$$\|\sqrt\Phi f\|_{H^1}\approx \left( \int \Phi (\vert f\vert^2+\vert \vN f\vert^2)\, dx\right)^{1/2}$$
b) $\Phi f\in W^{1,r}$ if and only if $f\in L^r(\Phi^r\, dx)$ and $\vN f\in L^r(\Phi^r\, dx)$; moreover we have
$$\| \Phi f\|_{W^{1,r}}\approx \left( \int \Phi^r (\vert f\vert^r+\vert \vN f\vert^r)\, dx\right)^{1/r}$$
\end{Lemma}

\pv{}
This is obvious since $\vert \vN \Phi\vert\leq C_1 \Phi^{3/2}\leq C_1\Phi$ and $\vert \vN(\sqrt\Phi)\vert =\frac 1 2 \frac{\vert \vN\Phi\vert}{\Phi} \sqrt\Phi\leq \frac 1 2 C_1\sqrt\Phi$. \Endproof{}

\begin{Lemma}\label{lem2}
If $\Phi \in \mathcal{A}_s$ then we have for all $\theta \in (0,1)$, $\Phi^\theta \in \mathcal{A}_p$ with $ \theta = \frac{p-1}{s-1}$. In particular, if a weight $\Phi$ satisfies $ (\textbf{H}3) $, we obtain $ \Phi \in \mathcal{A}_p$ with $ p = 1 + \frac{r-1}{r} =2 - \frac{1}{r} <2 $, and so $ \Phi \in \mathcal{A}_2$.
\end{Lemma}

\pv{}
As $ \frac{1}{p} =  \frac{1}{s}+ \frac{s-p}{ps} $, we find by the H\"older inequality
\begin{equation*}
\begin{split}
    (\int_{Q} & \Phi^{\frac{p-1}{s-1}} \, dx)^{\frac{1}{p}}  (\int_Q \Phi^{- (\frac{p-1}{s-1}) (\frac{1}{p-1} )} dx  )^{1-\frac{1}{p}} \\ 
    & = (\int_{Q} ( \Phi^{\frac{1}{s} } (\Phi^{-\frac{1}{s-1}})^{  \frac{s-p}{ps}} )^p \, dx)^{\frac{1}{p}}  (\int_Q \Phi^{- (\frac{p-1}{s-1}) (\frac{1}{p-1} )} dx )^{1-\frac{1}{p}} \\
    & \leq (\int_Q  \Phi \, dx)^{\frac{1}{s}}  (\int_{Q} \Phi^{-\frac{1}{s-1}} \, dx )^{\frac{1}{p}-\frac{1}{s} +1-\frac{1}{p}}
\end{split}
\end{equation*}
\Endproof{} 

Let us recall that for  a weight $w\in\mathcal{A}_q$ ($1<q<+\infty$), the Riesz transforms and the Hardy--Littlewood maximal function are bounded on $L^q(w\, dx)$.  We thus have the following inequalities:

\begin{Lemma}\label{lem3} Let $\Phi$ satisfy $(\textbf{H}1)$, $(\textbf{H}2)$ and $(\textbf{H}3)$. Then:
\\ a) for $j=1,\dots, d$, the Riesz transforms $R_j$ satisfy that $\|\sqrt \Phi R_jf\|_2\leq C \|\sqrt\Phi f\|_2$ and  $\|\sqrt \Phi R_jf\|_{H^1}\leq C \|\sqrt\Phi f\|_{H^1}$;
\\ b)   for $j=1,\dots, d$, the Riesz transforms $R_j$ satisfy that $\|  \Phi R_jf\|_r\leq C \| \Phi f\|_r$ and  $\|  \Phi R_jf\|_{W^{1,r}}\leq C \| \Phi f\|_{W^{1,r}}$;
\\ c) if $\mathbb{P}$ is the Leray projection operator on divergence-free vector fields, then for a vector field $\vu$ we have $\|\sqrt\Phi \mathbb{P} \vu\|_2\leq C \|\sqrt\Phi \vu\|_2$ and $\|\sqrt\Phi \mathbb{P} \vu\|_{H^1}\leq C \|\sqrt\Phi \vu\|_{H^1}$;
\\ d) if $d \in \{2,3,4 \}$, then for a vector field $\vu$ we have $$\|\sqrt \Phi \, \vu\|_{H^1} \approx \|\sqrt \Phi\, \vu\|_2+\|\sqrt\Phi \vN\cdot\vu\|_2+\|\sqrt\Phi\vN\wedge\vu\|_2.$$
  e) Let $\theta_\epsilon (x)=\frac 1{\epsilon^d}\theta(\frac x \epsilon )$, where $\theta\in\mathcal{D}(\mathbb{R}^d)$, $\theta$ is non-negative and radially decreasing and $\int\theta\, dx=1$. Then we have $\| \sqrt\Phi\, (\theta_\epsilon*f)\|_2\leq C\|\sqrt\Phi\, f\|_2$ and  $\| \sqrt\Phi\, (\theta_\epsilon*f)\|_{H^1}\leq C ( \|\sqrt\Phi\ f \|_{L^2} + \|\sqrt\Phi\ \vN f \|_{L^2})$ (where the constant $C$ does not depend on $\epsilon$ nor $f$).
\end{Lemma}

\pv{} a) is a consequence of $\Phi\in \mathcal{A}_2$ and of Lemma \ref{lem1} (since $\partial_k(R_jf)=R_j(\partial_kf)$). Similarly, b) is a consequence of $\Phi^r\in \mathcal{A}_r$ and of Lemma \ref{lem1}.

c) is a consequence of a): if $\vv=\mathbb{P}\vu$, then $v_j=\sum_{k=1}^d R_jR_k(u_k)$.

d) is a consequence of a): if $\mathcal{R}=(R_1,\dots, R_d)$, we have the identity $$-\Delta\vu=\vN\wedge(\vN\wedge\vu)-\vN (\vN\cdot \vu)$$ so that
$$ \partial_k\vu=R_k\mathcal{R}\wedge(\vN\wedge\vu)-R_k\mathcal{R}(\vN\cdot\vu).$$

e) is a consequence of $\Phi\in \mathcal{A}_2$ and of Lemma \ref{lem1}: Theorem 2.1.10 in Chapter 2 of \cite{Gr08} states that we have $\vert \theta_\epsilon*f\vert\leq \mathcal{M}_f$ (where $\mathcal{M}_f$ is the Hardy--Littlewood maximal function of $f$) and, similarly, $\vert  \partial_k(\theta_\epsilon*f)\vert\leq \mathcal{M}_{\partial_k f}$.\Endproof{}\\

A final lemma states that $\Phi$ is slowly decaying at infinity:

\begin{Lemma}\label{lem4} Let $\Phi$ satisfy $(\textbf{H}1)$ and $(\textbf{H}2)$. Then there exists a constant $C_3$ such that $$ \frac 1{(1+\vert x\vert)^2}\leq C_3 \Phi.$$

If $d=3$ and  $\Phi $ depends only on $r=\sqrt{x_1^2+x_2^2}$, then $$ \frac 1{(1+\vert r\vert)^2}\leq C_3 \Phi.$$
\end{Lemma}

\pv{} We define $x_0=\frac 1{\vert x\vert} x$ and $g(\lambda)=\Phi(\lambda x_0)$. We have
$$ g'(\lambda) =x_0\cdot \vN \Phi(\lambda x_0)\geq - C_1 (\Phi(\lambda x_0))^{3/2}=-C_1 g(\lambda)^{3/2}.$$
Thus
$$C_1 \lambda\geq -\int_0^\lambda g'(\mu)g(\mu)^{-3/2}\, d\mu= 2(g(\lambda)^{-1/2}-g(0)^{-1/2})$$ and we get
$$\Phi(x)^{-1/2} \leq \Phi(0)+\frac{C_1}2 \vert x\vert\leq \sqrt{C_3}(1+\vert x\vert).$$
If $\Phi$ depends only on $r$, we find that
$$ \frac 1{(1+\vert r\vert)^2}\leq C_3 \Phi(x_1,x_2,0)=C_3 \Phi(x).$$ \Endproof{}

\section{Proof of Theorem \ref{Thm1} (the case of $L^2(\Phi\, dx)$)}

\subsection{A priori controls}

Let $\phi \in \mathcal{D}(\mathbb{R}^d)$ be a real-valued test function which is equal to 1 in a neighborhood of $0$ and let $\phi_\epsilon(x)= \phi(\epsilon x)$. Let 
\begin{equation*}
    \vu_{0, \epsilon} = \mathbb{P}( \phi_\epsilon \vu_0 ) .
\end{equation*}
Thus, $\vu_{0, \epsilon}$ is divergence free and converges to $\vu_{0}$ in $L^2(\Phi \, dx)$
since $\Phi \in \mathcal{A}_2$.

Let $\theta_\epsilon (x)=\frac 1{\epsilon^d}\theta(\frac x \epsilon )$, where $\theta\in\mathcal{D}(\mathbb{R}^d)$, $\theta$ is non-negative and radially decreasing and $\int\theta\, dx=1$.
We denote $\vb_\epsilon = \vu_\epsilon*\theta_\epsilon$. Let $\vu_\epsilon$ be the unique global solution of the problem 
    \begin{equation*}  (NS_\epsilon) 
        \left\{ \begin{matrix} 
            \partial_t \vu_\epsilon  = \Delta \vu_\epsilon  -( \vb_\epsilon \cdot \vN) \vu_\epsilon - \vN p_\epsilon  
            \cr \cr \vN \cdot \vu_\epsilon =0, \phantom{space space} \vu_\epsilon(0,.)=\vu_{0, \epsilon}
        \end{matrix}\right.
    \end{equation*}
    which belongs to $\mathcal{C}([0,+\infty), L^2(\mathbb{R}^d))\cap L^2((0,+\infty),\dot H^1(\mathbb{R}^d))$.

We want to demonstrate that
\begin{equation}
\label{pc}
    \begin{split}
        \| \sqrt{\Phi}  \vu_\epsilon (t) \|_{L^2}^2 & + \int_0^t \| \sqrt{\Phi}  \vN \otimes  \vu_\epsilon  \|_{L^2}^2 \, ds \\
        & \leq \| \sqrt{\Phi} \vu_{0,\epsilon} \|_{L^2}^2 +  C_\Phi \int_0^t \| \sqrt{\Phi}  \vu_\epsilon  \|_{L^2}^2 + \| \sqrt{\Phi}  \vu_\epsilon  \|_{L^2}^{2d} \, ds , 
    \end{split}
\end{equation}
where $C_\Phi$ does not depend on $\epsilon$ nor on $\vu_0$. (When $d=4$, the inequality will hold only if $\| \sqrt{\Phi}  \vu_\epsilon (t) \|_{L^2}$ remains small enough).

Since   $\sqrt{\Phi},\vN \sqrt{\Phi} \in L^\infty$, pointwise multiplication by $\sqrt \Phi$ maps boundedly  $H^1$ to $H^1$ and $H^{-1}$ to $H^{-1}$. Thus,  $\sqrt{\Phi} \vu_\epsilon \in L^2 H^1$  and $\sqrt{\Phi} \partial_t \vu_\epsilon \in L^2 H^{-1}$, we can calculate $ \int \partial_t \vu_\epsilon \cdot \vu_\epsilon \Phi \,dx $ and obtain:
\begin{equation}
\label{eb}
    \begin{split}
        \int \frac{ |    \vu_\epsilon  (t,x) | ^2}2  \Phi \, dx & + \int_0^t  \int   | \vN  \otimes\vu_\epsilon   | ^2\, \,    \Phi dx\, ds \\
        =& \int \frac{ |   \vu_{0, \epsilon} (x) | ^2}2    \Phi\, dx 
        - \int_0^t  \int  (\vN  \otimes\vu_\epsilon) \, \cdot ( \vN  \Phi \otimes \vu_\epsilon) \, dx\, ds \\ & +  \int_0^t \int  ( \frac{| \vu_\epsilon |^2 }{2} \vb_\epsilon + p_\epsilon \vu_\epsilon)   \,  \cdot \vN   \Phi \, dx\, ds.
    \end{split}
\end{equation}

We use the fact that
$
    |\vN \Phi|  \leq C_0 \Phi^{\frac{3}{2}} \leq C_0 \Phi,
$
in order to control the following term
\begin{equation*}
    \left| - \int_0^t  \int ( \vN \otimes \vu_\epsilon) \, \cdot ( \vN  \Phi \otimes \vu_\epsilon)  dx\, ds \right| \leq \frac{1}{8} \int_0^t \|\sqrt\Phi\,  \vN\otimes \vu_\epsilon \|_{L^2}^2 + C \int_0^t \| \sqrt \Phi\, \vu_\epsilon \|_{L^2}^2. 
\end{equation*}

Now, we analyze the integrals containing the pressure term.  We distinguish two cases:

\begin{itemize}
    
    \item \textbf{Case 1:} $d=2$ and $r \in (1, 2]$, or $d=3$ and $r \in [\frac{6}{5}, 2]$, or $d=4$ and $r\in [\frac 4 3,2)$. For those values of $d$ and $r$ we have
  $$  0 \leq \frac{d}{2}-\frac{d}{2r} \leq 1 \text{ and }  \dot H^{\frac d 2-\frac d {2r}}\subset L^{2r}$$ and 
      $$  0 \leq \frac{d}{r}-\frac{d}{2} \leq 1\text{ and }  \dot H^{\frac d r-\frac d {2}}\subset L^{\frac r {r-1}}. $$

    Using the continuity of the Riesz transforms on $L^r(\Phi^r dx)$, 
    \begin{equation*}
    \begin{split}
        \int_0^t   \int  & ( \frac{| \vu_\epsilon |^2 |\vb_\epsilon| }{2} + |p_\epsilon| |\vu_\epsilon| )   \,  | \vN   \Phi | \, dx\, ds  \leq \int_0^t \| \Phi ( | \vu_\epsilon |\,  | \vb_\epsilon | + |p_\epsilon| ) \|_r  \| \sqrt{\Phi} \vu_\epsilon \|_{\frac{r}{r-1}}  \\
        & \leq  C  \int_0^t \| \sqrt{\Phi}  \vu_\epsilon \|_{2r} \| \sqrt{\Phi}   \vb_\epsilon  \|_{2r}  \| \sqrt{\Phi} \vu_\epsilon \|_{\frac{r}{r-1}} ds
    \end{split}
    \end{equation*}
    Using the Sobolev embedding $\dot H^{\frac d 2-\frac d {2r}}\subset L^{2r}$, the fact that $ |\vN \sqrt{\Phi}| \leq C \sqrt{\Phi} $, and   the continuity of the maximal function operator on $L^2 (\Phi dx)$, we have
    \begin{equation*}
    \begin{split}
        \| & \sqrt{\Phi}  \vb_\epsilon  \|_{2r} \\
        & \leq C \| \sqrt{\Phi}  \vb_\epsilon  \|_{2}^{1-(\frac{d}{2}- \frac{d}{2r})} \| \vN \otimes(\sqrt{\Phi}  \vb_\epsilon  ) \|_{2}^{\frac{d}{2}- \frac{d}{2r}} \\
        & \leq C' \| \sqrt{\Phi}  \vb_\epsilon  \|_{2}^{1-(\frac{d}{2}- \frac{d}{2r})} (  \| \sqrt{\Phi} \vb_\epsilon  \|_{2} + \| \sqrt{\Phi}  \vN \otimes  \vb_\epsilon  \|_{2})^{\frac{d}{2}- \frac{d}{2r}}
         \\
        & \leq C'' \| \sqrt{\Phi}  \vu_\epsilon  \|_{2}^{1-(\frac{d}{2}- \frac{d}{2r})} ( \| \sqrt{\Phi} \vu_\epsilon  \|_{2} + \| \sqrt{\Phi}  \vN \otimes\vu_\epsilon  \|_{2})^{\frac{d}{2}- \frac{d}{2r}} ,
    \end{split}
    \end{equation*} and 
     \begin{equation*}
        \|  \sqrt{\Phi}  \vu_\epsilon  \|_{2r} \\
      \leq C \| \sqrt{\Phi}  \vu_\epsilon  \|_{2}^{1-(\frac{d}{2}- \frac{d}{2r})} ( \| \sqrt{\Phi} \vu_\epsilon  \|_{2} +  \| \sqrt{\Phi}  \vN \otimes\vu_\epsilon  \|_{2})^{\frac{d}{2}- \frac{d}{2r}}.
  \end{equation*}

 Using the embedding $ \dot H^{\frac d r-\frac d {2}}\subset L^{\frac r {r-1}}$, we also have
    \begin{equation*}
    \begin{split}
        \| & \sqrt{\Phi} \vu_\epsilon \|_{\frac{r}{r-1}}  \\
        & \leq  C   \| \sqrt{\Phi}  \vu_\epsilon \|_{2}^{1-( \frac{d}{r}-\frac{d}{2}  ) } \| \vN \otimes( \sqrt{\Phi}  \vu_\epsilon ) \|_{L^2}^{\frac{d}{r}-\frac{d}{2}  }   \\
        & \leq  C   \| \sqrt{\Phi}  \vu_\epsilon \|_{2}^{1-( \frac{d}{r}-\frac{d}{2}  ) } (\| \sqrt{\Phi}  \vu_\epsilon \|_{2} +  \| \sqrt{\Phi}  \vN \otimes  \vu_\epsilon \|_{L^2})^{ \frac{d}{r}-\frac{d}{2}  } . 
    \end{split}
    \end{equation*}
    
    Hence, we find 
     \begin{equation*}
    \begin{split}
        \int_0^t   \int  & ( \frac{| \vu_\epsilon |^2 |\vb_\epsilon| }{2} + |p_\epsilon| |\vu_\epsilon| )   \,  | \vN   \Phi | \, dx \, ds \\
        & \leq  C \int_0^t   \| \sqrt{\Phi}  \vu_\epsilon \|_{2}^{3- \frac{d}{2} } (\| \sqrt{\Phi}  \vu_\epsilon \|_{2} +  \| \sqrt{\Phi}  \vN \otimes  \vu_\epsilon \|_{L^2})^{ \frac{d}{2}  }\, ds . 
          \end{split}
    \end{equation*}

Using the Young inequality, we then find for $d=2$ or $d=3$
    \begin{equation*}
    \begin{split}
        \int_0^t   \int  & ( \frac{| \vu_\epsilon |^2 |\vb_\epsilon| }{2} + |p_\epsilon| |\vu_\epsilon| )   \,  | \vN   \Phi | \, dx \, ds \\
        & \leq \frac{1}{8}\int_0^t \| \sqrt{\Phi}\vN\otimes  \vu_\epsilon  \|_{L^2}^2 \, ds + C_\Phi \int_0^t \| \sqrt{\Phi} \vu_\epsilon \|_{L^2}^2 +  \| \sqrt{\Phi} \vu_\epsilon \|_{L^2}^{\frac{12-2d}{4-d}} \, ds,
    \end{split}
    \end{equation*}
    where, as $d \in \{2,3\}$, we have $\frac{12-2d}{4-d}= 2d$.
    
    When $d=4$, provided that $\|\sqrt\Phi\, \vu_\epsilon\|_2<\epsilon_0$ with $C\epsilon_0<\frac 1 8$
    we find 
      \begin{equation*}
    \begin{split}
        \int_0^t   \int  & ( \frac{| \vu_\epsilon |^2 |\vb_\epsilon| }{2} + |p_\epsilon| |\vu_\epsilon| )   \,  | \vN   \Phi | \, dx \, ds \\
        & \leq \frac{1}{8} \int_0^t \| \sqrt{\Phi} \vN\otimes \vu_\epsilon  \|_{L^2}^2  \, ds+  \frac{1}{8}\int_0^t \| \sqrt{\Phi} \vu_\epsilon \|_{L^2}^2  \, ds,
    \end{split}
    \end{equation*}
     
    \item \textbf{Case 2:} $d=3$ and $r \in ( 1, \frac{6}{5} )$, or $d=4$ and $r\in (1,\frac 4 3)$. Let $q = \frac{dr}{d-r}$; for those values of $d$, $r$ and $q$,  we have
  $$  W^{1,r}\subset L^{q}$$
  $$  0 \leq d-\frac{d}{r} \leq 1\text{ and }  \dot H^{  d (1-\frac 1 {r})}\subset L^{ \frac{2r}{2-r}}. $$ and 
      $$  0 \leq \frac{d}{r}-\frac{d}{2}-1 \leq 1\text{ and }  \dot H^{\frac d r-\frac d {2}-1}\subset L^{\frac q{q-1}}. $$
    
      Using the continuity of the Riesz transforms on $L^r(\Phi^r dx)$, we have

    \begin{equation*}
    \begin{split}
        & \int_0^t  \int  ( \frac{| \vu_\epsilon |^2 |\vb_\epsilon  | }{2} + |p_\epsilon| |\vu_\epsilon | )   \,  | \vN   \Phi | \, dx\, ds \\
        \leq & \int_0^t \| \Phi  | \vu_\epsilon|^2  \|_{q}  \| \sqrt{\Phi} \vb_\epsilon  \|_{\frac{q}{ q-1}} ds  + \int_0^t  \| \Phi p_\epsilon  \|_q  \| \sqrt{\Phi} \vu_\epsilon \|_{\frac{q}{q-1}} ds \\
        \leq  & C \int_0^t \| \Phi  | \vu_\epsilon |^2  \|_{W^{1,r}}  \| \sqrt{\Phi} \vb_\epsilon \|_{\frac{q}{q-1}} ds + \sum_{ij} \int_0^t  \| \Phi b_{\epsilon, i} u_{\epsilon, j}   \|_{W^{1,r}}  \| \sqrt{\Phi} \vu_\epsilon \|_{\frac{q}{q-1}} ds.
        \end{split}
        \end{equation*} 
         We have 
        \begin{equation*}
        \begin{split}
        & \|  \Phi b_{\epsilon, i}  u_{\epsilon, j}   \|_{W^{1,r}} \\
        \leq &  \|\Phi  b_{\epsilon, i}  u_{\epsilon, j}   \|_r+\sum_{k} ( \| b_{\epsilon,i}  u_{\epsilon,j} \, \partial_k \Phi \|_{L^r} + \| \Phi  \,  b_{\epsilon,i} \, \partial_k  u_{\epsilon,j}   \|_{L^r} + \| \Phi  \,  u_{\epsilon,i} \, \partial_k  b_{\epsilon,j} \|_{L^r})  \\
        \leq &  C  ( \| \sqrt{\Phi}  \vu_\epsilon \|_{\frac{2r}{2-r}}   \| \sqrt{\Phi} \vb_\epsilon  \|_{2}+  \| \sqrt{\Phi} \vb_\epsilon  \|_{\frac{2r}{2-r}} \| \sqrt{\Phi}  \vN \otimes\vu_\epsilon \|_{2} + \| \sqrt{\Phi} \vu_\epsilon \|_{\frac{2r}{2-r}}    \| \sqrt{\Phi} \vN\otimes\vb_\epsilon  \|_2),
        \\ 
        \leq & C'  ( \| \sqrt{\Phi}  \vu_\epsilon \|_{L^2} + \| \sqrt{\Phi}  \vN \otimes\vu_\epsilon \|_{L^2} ) (\| \sqrt{\Phi} \vu_\epsilon \|_{ \dot H^{  d (1-\frac 1 {r})}} + \| \sqrt{\Phi} \vb_\epsilon  \|_{\dot H^{  d (1-\frac 1 {r})}}).
    \end{split}
    \end{equation*}
We have 
    \begin{equation*}
    \begin{split}
        \| & \sqrt{\Phi}  \vb_\epsilon  \|_{ \dot H^{  d (1-\frac 1 {r})}}\\
        & \leq C \| \sqrt{\Phi}  \vb_\epsilon  \|_{2}^{1-(d- \frac{d}{r})} \| \vN\otimes (\sqrt{\Phi}  \vb_\epsilon  ) \|_{2}^{d- \frac{d}{r}} \\
        & \leq C' \| \sqrt{\Phi}  \vb_\epsilon  \|_{2}^{1-(d- \frac{d}{r})} (   \| \sqrt{\Phi} \vb_\epsilon  \|_{2} + \| \sqrt{\Phi}  \vN  \otimes\vb_\epsilon  \|_{2})^{d- \frac{d}{r}}
         \\
        & \leq C'' \| \sqrt{\Phi}  \vu_\epsilon  \|_{2}^{1-(d- \frac{d}{r})} ( \| \sqrt{\Phi} \vu_\epsilon  \|_{L^2} + \| \sqrt{\Phi} \vN \otimes \vu_\epsilon  \|_{L^2})^{d- \frac{d}{r}} ,
    \end{split}
    \end{equation*}
    and finally we get
   \begin{equation*}
    \begin{split}
  \sum_{i,j} \| & \Phi b_{\epsilon, i}  u_{\epsilon, j}   \|_{W^{1,r}}  +  \|  \Phi  \vert u_{\epsilon}\vert^2 \|_{W^{1,r}}    \\
        & \leq C \| \sqrt{\Phi}  \vu_\epsilon  \|_{2}^{1-(d- \frac{d}{r})} ( \| \sqrt{\Phi} \vu_\epsilon  \|_{L^2} + \| \sqrt{\Phi} \vN \otimes \vu_\epsilon  \|_{L^2})^{1+d- \frac{d}{r}} .
      \end{split}
    \end{equation*}
    
    On the other hand, we have
      \begin{equation*}
    \begin{split}
        \| & \sqrt{\Phi}  \vb_\epsilon  \|_{\frac{q}{q-1}} \\
        & \leq C \| \sqrt{\Phi}  \vb_\epsilon  \|_{2}^{2-(\frac d r- \frac{d}{2})} \| \vN \otimes(\sqrt{\Phi}  \vb_\epsilon  ) \|_{2}^{\frac d r- \frac{d}{2}-1} \\
        & \leq C'   \| \sqrt{\Phi}  \vu_\epsilon  \|_{2}^{2-(\frac d r- \frac{d}{2})}   ( \| \sqrt{\Phi} \vu_\epsilon  \|_{L^2} + \| \sqrt{\Phi} \vN\otimes  \vu_\epsilon  \|_{L^2})^{\frac d r- \frac{d}{2}-1}.
    \end{split}
    \end{equation*}
    
 Hence, we find again
     \begin{equation*}
    \begin{split}
        \int_0^t   \int  & ( \frac{| \vu_\epsilon |^2 |\vb_\epsilon| }{2} + |p_\epsilon| |\vu_\epsilon| )   \,  | \vN   \Phi | \, dx \, ds \\
        & \leq  C \int_0^t   \| \sqrt{\Phi}  \vu_\epsilon \|_{2}^{3- \frac{d}{2} } (\| \sqrt{\Phi}  \vu_\epsilon \|_{2} +  \| \sqrt{\Phi}  \vN \otimes  \vu_\epsilon \|_{L^2})^{ \frac{d}{2}  }\, ds . 
          \end{split}
    \end{equation*} and we conclude in the same way as for the first case. 
\end{itemize}

In the Case 1 and Case 2, we have found
\begin{equation*}
    \begin{split}
        \int_0^t   \int  & ( \frac{| \vu_\epsilon |^2 |\vb_\epsilon| }{2} + |p_\epsilon| |\vu_\epsilon| )   \,  | \vN   \Phi | \, dx \, ds \\
        & \leq \frac{1}{8} \| \sqrt{\Phi}  \vu_\epsilon  \|_{L^2}^2  + C_\Phi \int_0^t \| \sqrt{\Phi} \vu_\epsilon \|_{L^2}^2 +  \| \sqrt{\Phi} \vu_\epsilon \|_{L^2}^{2d} \, ds.
    \end{split}
\end{equation*}

From these controls, we get inequality   \eqref{eb}, and thus inequality
  \eqref{pc}. Inequality   \eqref{pc} gives us a control on the size of $\|\sqrt\Phi\, \vu_\epsilon\|_2$ on an interval of time that does not depend on $\epsilon$:

\begin{Lemma}
\label{lem5} 
If $\alpha$ is a continuous non-negative function on $[0,T)$ which satisfies, for three constants $A,B \in (0, + \infty)$  and $ b \in [1, \infty) $,
$$ \alpha(t)\leq A + B\int_0^t \alpha(s) + \alpha(s)^b\, ds.$$ 
Let $0<T_1<T$ and $T_0=\min( T_1, \frac 1{3^b (A ^{b-1}+ (BT_1)^{b-1})})$. We have, for every $t\in [0,T_0]$,  $\alpha(t)\leq  3A$.
\end{Lemma}

\pv     We try to estimate the first time $T^*<T_1$ (if it exists)  for which we have 
$$\alpha(T^*) =3A.$$ We have 
$$ \alpha\leq \frac A{BT_1} + (\frac   {BT_1}A)^{b-1} \alpha^b.$$ 
We thus find
$$ \alpha(T^*)\leq 2A + T^* (3A)^b (1+(\frac  {BT_1}A)^{b-1} )$$
and thus
$$ T^* 3^b (A^{b-1}+ (BT_1)^{b-1})\geq 1.$$ \Endproof\\

By Lemma \ref{lem5} and \eqref{pc}, we thus  find that  there exists a constant $C_\Phi \geq 1$ such that if $T_0 $ satisfies
\begin{itemize}
\item[$\bullet$] if $d=2$, $ C_\Phi \left(1+\|\vu_0\|_{L^2(\Phi dx)}^2  \right) \, T_0 \leq 1$
\item[$\bullet$] if $d=3$, $ C_\Phi \left(1+\|\vu_0\|_{L^2(\Phi dx)}^2  \right)^2 \, T_0 \leq 1$
\item[$\bullet$] if $d=4 $ and $\|\vu_0\|_{L^2(\Phi\, dx)}\leq \frac 1{C_\Phi}$, $ C_\Phi  \, T_0 \leq 1$
\end{itemize} 
then
\begin{equation}\label{ineq_energy}
        \sup_{0\leq t\leq T_0} \|\ \vu_\epsilon(t,.)\|_{L^2(\Phi dx)}^2 +  \int_0^{T_0} \|\vN \otimes\vu_\epsilon \|_{L^2(\Phi \, dx)}^2\, ds  \leq
        C_\Phi (1 + \|\vu_0\|_{L^2(\Phi\,  dx)}^2 ).
\end{equation}


\subsection{Passage to the limit and local existence}
\label{plle}

We know 
that $\vu_\epsilon$ is bounded in $L^\infty((0,T_0 ), L^2(\Phi \, dx))$  and $ \vN \otimes\vu_\epsilon$ is bounded  in $L^2((0,T_0 ), L^2(\Phi \, dx))$. This will alow us to use a simple variant of the Aubin--Lions theorem:

\begin{Lemma}[Aubin--Lions theorem] \label{lem6} Let $s >0$, $1<q$ and $\sigma <0$. Let $ (f_n)$ be a sequence of functions on $(0,T)\times \mathbb{R}^d$ such that, for all $T_0\in (0,T)$ and all $\varphi\in\mathcal{D}(\mathbb{R}^d)$,
  \begin{itemize}
  \item[$\bullet$]  $\varphi f_n$ is bounded in $L^2((0,T_0), H^s)$ 
  \item[$\bullet$]   $\varphi \partial_t f_n$ is bounded in $L^q((0,T_0), H^\sigma)$ .
\end{itemize}

Then, there exists a subsequence $(f_{n_k})$ such that $f_{n_k}$ is strongly convergent in $L^2_{\rm loc}([0,T)\times \mathbb{R}^d)$. More precisely: if we denote $f_\infty$  the limit, then for all $T_0\in (0,T)$ and all $R_0>0$,
$$\lim_{n_k\rightarrow +\infty} \int_0^{T_0} \int_{ |  x | \leq R_0} | f_{n_k}-f_\infty |^2 \, dx\, dt=0.$$
\end{Lemma}
 For a proof of the Lemma, see  \cite{BF12, LR16}.\\ 

We want to verify that $\varphi\partial_t\vu_\epsilon$ is bounded in $ L^\alpha((0,T_0 ), H^{-s})$ for some $s \in (-\infty,0)$ and some $\alpha>1$.

In   Case 1, we have  that $\Phi \vb_\epsilon\otimes\vu_\epsilon$  and $\Phi p_\epsilon =\sum_{i=1}^3\sum_{j=1}^3 R_iR_j(b_{\epsilon,i}u_{\epsilon,j})  $
are bounded in $L^{\alpha_1}((0,T_0 ),L^{r})$, where $\alpha_1 = \frac{2r}{ dr-d}$, so that $\alpha_1 \in [2, \infty)$ if $d=2$, $\alpha_1 \in [\frac{4}{3}, 4]$ if $d=3$  and $\alpha_1 \in (1,2]$ if $d=4$.

In   Case 2,  we have  that $\Phi \vb_\epsilon\otimes\vu_\epsilon$  and  $\Phi p_\epsilon $ are bounded in $L^{\alpha_2}((0,T_0 ),W^{1,r})$, where $ \alpha_2 = \frac{2r}{r + dr-d} $  and thus it is bounded in $L^{\alpha_2} L^{q}$, with $q = \frac{dr}{d-r}$.  We have  $\alpha_2 \in (\frac{4}{3}, 2)$ if $d=3$  and $\alpha_2 \in (1,2)$ if $d=4$.

Let $\varphi\in \mathcal{D}(\mathbb{R}^d)$. We have that $\varphi \vu_\epsilon $ is bounded in $L^2((0,T_0 ), H^1)$;   
moreover, writing
\begin{equation*}
    \begin{split}
        \partial_t \vu_\epsilon &= \Delta \vu_\epsilon  - \left( \sum_{j=1}^3 \partial_j(b_{\epsilon,j}\vu_\epsilon) +\vN p_{\epsilon}\right)
    \end{split}
\end{equation*} and using the embeddings $L^r\subset  \dot H^{\frac d {2} - \frac d r} \subset H^{-1}$ (in Case 1) or $L^{\frac{dr}{d-r}}\subset H^{-(\frac d r-\frac d {2}-1)}\ \subset H^{-1}$ (in Case 2)
we see that $\varphi \partial_t \vu_\epsilon$ is bounded in $L^{\alpha_i}((0,T_0), H^{-2})$.
 
Thus, by the Aubin--Lions theorem, there exist  $\vu$  and a  sequence $(\epsilon_k)_{k\in\mathbb{N}}$ converging to $0$ such that $\vu_{\epsilon_k}$ converges strongly  to $\vu$ in  $L^2_{\rm loc} ([0,T_0 )\times\mathbb{R}^3)$: for every  $\tilde  T \in (0,T_0 )$ and every $R>0$, we have
\begin{equation*}
        \lim_{k\rightarrow +\infty} \int_0^{\tilde T} \int_{ |  y | <R}  |  \vu_{\epsilon_k}-\vu  |^2\, dx\, ds  =0   .
\end{equation*}

Then, we have that $\vu_{\epsilon_k}$  converge *-weakly to $\vu $ in $L^\infty((0,T_0 ), L^2(\Phi dx))$, $\vN \otimes\vu_{\epsilon_k}$   converges weakly to $\vN\otimes\vu$  in $L^2((0,T_0 ),L^2(\Phi dx))$, and $ \vu_{\epsilon_k}$  converges weakly to $ \vu$ in $L^3((0,T_0 ), L^3(\Phi^{\frac{3}{2}} dx))$.
We deduce that $ \vb_{\epsilon_k}\otimes \vu_{\epsilon_k}$  
is weakly convergent in $(L^{6/5}L^{6/5})_{\rm loc}$ to $\vb\otimes\vu$ and thus in $\mathcal{D}'((0,T_0 )\times \mathbb{R}^d)$; as in   Case 1, it is bounded in $L^{\alpha_1}((0,T_0 ),L^{r})$, and in   Case 2 it is bounded in $L^{\alpha_2}((0,T_0 ),W^{1,r})$, it is weakly convergent in these spaces respectively (as $\mathcal{D}$ is dense in their dual spaces).

By the continuity of the  Riesz transforms on $L^{r}(\Phi^{r} dx)$ and on $W^{1,r}(\Phi^{r} dx)$, we find that in the Case 1 and Case 2, $  p_{\epsilon_k}$ is convergent to the distribution $p= \sum_{i=1}^3\sum_{j=1}^3 R_iR_j(u_{i}u_{j})$. We have obtained
\begin{equation*}
    \partial_t \vu = \Delta \vu + ( \vu \cdot \vN)\vu - \vN p.
\end{equation*} 

Moreover, we have seen that $\partial_t \vu$  is locally in $L^1 H^{-2}$, and thus $\vu$ has representative such that  $t\mapsto \vu(t,.)$
is continuous from $[0,T_0 )$ to $\mathcal{D}'(\mathbb{R}^d)$ and coincides with $\vu(0,.)+\int_0^t \partial_t \vu\, ds$.

In the sense of distributions, we have $$\vu(0,.)+\int_0^t \partial_t \vu\, ds=\vu=\lim_{k\rightarrow +\infty} \vu_{\epsilon_k}=\lim_{k\rightarrow +\infty}
\vu_{0,\epsilon_k}+ \int_0^t \partial_t \vu_{n_k}\, ds=\vu_{0}+\int_0^t \partial_t\vu\, ds,$$
hence, $\vu(0,.)=\vu_{0}$, and $\vu$ is a solution of $(NS)$.

Now, we want to prove the energy balance.  In the case of dimension 2, we remark that, since $\sqrt{\Phi}\vu \in L^\infty L^2 \cap L^2 H^1$, we have by interpolation that $\sqrt{\Phi}\vu \in L^4 L^4$, and then we can define
$ ((\vu \cdot \vN) \vu) \cdot \vu$. The equality \begin{equation*}
    \partial_t(\frac { | \vu | ^2}2)=\Delta(\frac { | \vu | ^2}2)- | \vN \vu| ^2- \vN\cdot\left( \frac{ | \vu | ^2}2\vu \right)-\vN\cdot(p\vu) 
\end{equation*} is then easy to prove.

Let us consider the case $d = 3$. We define 
\begin{equation*}
    A_\epsilon=  - \partial_t(\frac { | \vu_\epsilon  | ^2}2)+\Delta(\frac { | \vu_\epsilon  | ^2}2)-\vN\cdot\left( \frac{ | \vu_\epsilon  | ^2}2\vu_\epsilon \right)-\vN\cdot(p_\epsilon \vu_\epsilon ) = | \vN \otimes\vu_\epsilon  | ^2.
\end{equation*}
As $u_{\epsilon_k}$ is locally strongly convergent in $L^2 L^2$; and locally bounded in $L
^\infty L^2$, it is then locally strongly convergent in $L^{p'} L^2$, with $p' < \infty$. Then, as $\sqrt{\Phi} \vN\otimes\vu_\epsilon$ is bounded  in $L^2((0,T), L^2)$, by the Gagliardo-Nirenberg interpolation inequalities we obtain $\vu_{\epsilon_k}$ is locally strongly convergent in $L^{p'} L^{q'}$ with $ \frac{2}{p'}+ \frac{3}{q'} > \frac{d}{2} $. 

In   Case 1, we know that $p_{\epsilon_k}$ is locally weakly convergent in $L^{\alpha_1} L^r$ and by the remark above, $\vu_{\epsilon_k}$ is locally strongly convergent in $L^{\frac{\alpha_1}{\alpha_1 -1}} L^{\frac{r}{r-1}}$, and hence $p_{\epsilon_k} \vu_{\epsilon_k}$ converges in the sense of distributions.

In   Case 2, we know that $p_{\epsilon_k}$ is locally weakly convergent in $L^{\alpha_2} L^{q}$ and and by the remark above, $\vu_{\epsilon_k}$ is locally strongly convergent in $L^{\frac{\alpha_2}{\alpha_2-1}} L^{\frac{q}{q -1}}$, and hence $p_{\epsilon_k} \vu_{\epsilon_k}$ converges in the sense of distributions.

Thus, $A_{\epsilon_k}$ is convergent in $\mathcal{D}'((0,T)\times\mathbb{R}^3)$ to
\begin{equation*}
    A=  - \partial_t(\frac { | \vu | ^2}2)+\Delta(\frac { | \vu | ^2}2)-\vN\cdot\left( \frac{ | \vu | ^2}2\vu \right)-\vN\cdot(p \vu ),
\end{equation*}
 and $A =\lim_{k \rightarrow +\infty}  | \vN \otimes\vu_{\epsilon_k} | ^2  $. If $\theta\in \mathcal{D}((0,T)\times \mathbb{R}^d)$ is non-negative, we have that $\sqrt\theta \vN\otimes \vu_{\epsilon_k}$ is weakly convergent in $L^2L^2$ to $\sqrt\theta \vN\otimes \vu$, so that
 \begin{equation*}
    \iint A \theta\, dx\, ds=\lim_{\epsilon_k\rightarrow +\infty}\iint A_{\epsilon_k} \theta\, dx\, ds= \lim_{k\rightarrow +\infty}\iint    | \vN\otimes \vu_{\epsilon_k} | ^2 \theta\, dx\, ds \geq  \iint    | \vN\otimes \vu | ^2 \theta\, dx\, ds .
 \end{equation*}
Hence, there   exists a non-negative locally finite measure $\mu$ on $(0,T)\times\mathbb{R}^3$ such that $A= | \vN \vu | ^2 +\mu$, i.e. such that
\begin{equation*}
    \partial_t(\frac { | \vu | ^2}2)=\Delta(\frac { | \vu | ^2}2)- | \vN \vu| ^2- \vN\cdot\left( \frac{ | \vu | ^2}2\vu \right)-\vN\cdot(p\vu) -\mu.
\end{equation*}


\subsection{Strong  convergence to the initial data}

We use again inequalities \eqref{pc} and \eqref{ineq_energy}. We know that on $(0,T_0)$ we have a control of $\|\vu_\epsilon\|_{L^2(\Phi \, dx)}$ that holds uniformly in $\epsilon$ and $t$. Thus, inequality \eqref{pc} gives us
$$ \|\vu_\epsilon(t,.)\|_{L^2(\Phi \, dx)}\leq \|\vu_{0,\epsilon}\|_{L^2(\Phi \, dx)} + C_\Phi t (1+\|\vu_0\|_{L^2(\Phi \, dx)}^2+\|\vu_0\|_{L^2(\Phi \, dx)}^{2d}).$$

Since $ \vu_{\epsilon_k}= 
\vu_{0, \epsilon_k}+ \int_0^t \partial_t \vu_{\epsilon_k}\, ds $, we see that $\vu_{\epsilon_k}(t,.)$ is convergent to $\vu(t,.)$ in $\mathcal{D}'(\mathbb{R}^d)$, hence is weakly convergent in $L^2(\Phi\, dx)$ (as it is bounded in $L^2(\Phi dx)$); on the other hand, $\vu_{0,\epsilon_k}$ is  strongly convergent to $\vu_0$  in $L^2(\Phi\, dx)$. Thus, we have
$$ \|\vu(t,.)\|_{L^2(\Phi \, dx)}\leq \|\vu_{0}\|_{L^2(\Phi \, dx)} + C_\Phi t (1+\|\vu_0\|_{L^2(\Phi \, dx)}^2+\|\vu_0\|_{L^2(\Phi \, dx)}^{2d}).$$ In particular,
$$ \limsup_{t\rightarrow 0} \|\vu(t,.)\|_{L^2(\Phi \, dx)}\leq \|\vu_{0}\|_{L^2(\Phi \, dx)} .$$

Moreover, we have $ \vu= 
\vu_{0}+ \int_0^t \partial_t \vu\, ds $, so that $\vu (t,.)$ is convergent to $\vu_0$ in $\mathcal{D}'(\mathbb{R}^d)$, hence is weakly convergent in $L^2(\Phi\, dx)$. Thus, we have
$$\|\vu_{0}\|_{L^2(\Phi \, dx)}\leq  \liminf_{t\rightarrow 0} \|\vu(t,.)\|_{L^2(\Phi \, dx)}.$$

This gives $    \|\vu_{0} \|_{L^2(\Phi dx)}^2 =   \lim_{t\rightarrow 0}  \|\vu(t,. )\|_{L^2(\Phi dx)}^2 $, which allows to turn the weak convergence into a strong convergence.
\Endproof{}


\subsection{Global existence using a scaling argument}

Let $\lambda>0$, then $\vu_{\epsilon}$ is a solution of the Cauchy initial value problem for the approximated Navier--Stokes equations  $(NS_\epsilon)$ on $(0,T)$ with initial value $\vu_{0, \epsilon}$  if and only if $\vu_{{\epsilon}, \lambda} (t,x)=\lambda \vu_{\epsilon} (\lambda^2 t,\lambda x)$ is a solution  for the approximated Navier--Stokes equations  $(NS_{\lambda\epsilon})$  on $(0,T/\lambda^2)$ with initial value $\vu_{0,{\epsilon},\lambda}(x)=\lambda \vu_{0, {\epsilon}}(\lambda x)$. We shall write $\vu_{0,\lambda}=\lambda\vu_0(\lambda x)$.
 
We have seen that
\begin{equation*}
    \begin{split}
        \| \sqrt{\Phi}  \vu_{\epsilon , \lambda } (t) \|_{L^2}^2 & + \int_0^t \| \sqrt{\Phi}  \vN   \otimes\vu_{\epsilon , \lambda }  \|_{L^2}^2 \\
        & \leq \| \sqrt{\Phi} \vu_{0, \epsilon , \lambda } \|_{L^2}^2 +  C_\Phi \int_0^t \| \sqrt{\Phi}  \vu_{\epsilon , \lambda }  \|_{L^2}^2 + \| \sqrt{\Phi}  \vu_{\epsilon , \lambda }  \|_{L^2}^{2d}  \, ds
    \end{split}
\end{equation*}
(under the extra condition, when $d=4$, that $\| \sqrt{\Phi}  \vu_{\epsilon , \lambda } (t) \|_{L^2}$ remains smaller than $\epsilon_0$).

By Lemma \ref{lem5}, we thus  found that  there exists a constant $C_\Phi \geq 1$ such that if $T_\lambda $ satisfies
\begin{itemize}
\item[$\bullet$] if $d=2$, $ C_\Phi \left(1+\|\vu_{0,\lambda}\|_{L^2(\Phi dx)}^2  \right) \, T_\lambda = 1$
\item[$\bullet$] if $d=3$, $ C_\Phi \left(1+\|\vu_{0,\lambda}\|_{L^2(\Phi dx)}^2  \right)^2 \, T_\lambda = 1$
\item[$\bullet$] if $d=4 $ and $\|u_{0,\lambda}\|_{L^2(\Phi\, dx)}\leq \frac 1{C_\Phi}$, $ C_\Phi  \, T_\lambda = 1$
\end{itemize} 
then
\begin{equation}\label{ineq_energyter}
        \sup_{0\leq t\leq T_\lambda} \|\ \vu_{\epsilon,\lambda}(t,.)\|_{L^2(\Phi dx)}^2 +  \int_0^{T_\lambda} \|\vN \otimes\vu_{\epsilon,\lambda} \|_{L^2(\Phi \, dx)}^2\, ds  \leq
        C_\Phi (1 + \|\vu_{0,\lambda}\|_{L^2(\Phi\,  dx)}^2 ).
\end{equation}

It gives that the solutions $\vu_\epsilon$ are controlled, uniformly in $\epsilon$, on $(0, \lambda^2 T_\lambda)$ since for $t \in (0, T_\lambda)$,
$$ \int \vert \vu_{{\epsilon}, \lambda}(t,x)\vert^2 \Phi(x)\, dx =  \int \vert \vu_\epsilon(\lambda^2 t, y)\vert^2 \Phi(\frac y \lambda) \lambda^{2-d} \, dy\geq  \lambda^{2-d}\int \vert \vu_{{\epsilon}}(\lambda^2 t,x)\vert^2 \Phi(x)\, dx 
  $$
 and
 \begin{equation*} 
    \begin{split}   \int_0^{T_\lambda}  
    \int \vert \vN \otimes\vu_{\epsilon, \lambda  }(t,x) \vert^2 \Phi(x)\, dx\, dt=&    \int_0^{\lambda^2 T_\lambda}  
    \int \vert \vN \otimes\vu_{\epsilon, \lambda  }(s,y) \vert^2 \Phi(\frac y \lambda) \lambda^{2-d}\, dy\, ds
    \\ \geq&  \lambda^{2-d}  \int_0^{\lambda^2 T_\lambda}  
    \int \vert \vN \otimes\vu_{\epsilon  }(t,x) \vert^2 \Phi(  x )\, dx\, dt
    \end{split}
 \end{equation*} \begin{equation*} 
    \begin{split}   \int_0^{T_\lambda}  
    \int \vert \vN \otimes\vu_{\epsilon, \lambda  }(t,x) \vert^2 \Phi(x) \, dx\, dt \geq&  C_{\lambda} \int_0^{\lambda^2 T_\lambda}\| \vN\otimes \vu_{\epsilon} \|_{L^2(\Phi dx)}^2\, ds.
    \end{split}
 \end{equation*} 
 
Moreover, we have $ \lim_{ \lambda \rightarrow +\infty} \|\vu_{0,\lambda}\|_{L^2(\Phi\, dx)}=0$ when $d=4$ and  $\lim_{ \lambda \rightarrow +\infty} \lambda^2 T_\lambda = + \infty$ when $2\leq d\leq 4$. Indeed, we have
$$ \int  \lambda^2\vert \vu_0(\lambda x)\vert^2 \Phi(x)\, dx=\lambda^{2-d} \int   \vert \vu_0(  x)\vert^2 \Phi(\frac x \lambda)\, dx = \lambda^{4-d} \int   \vert \vu_0(  x)\vert^2 \frac{\Phi(\frac x \lambda)}{\lambda^2\Phi(x)} \Phi(x)\, dx $$
Since $\frac{\Phi(\frac x \lambda)}{\lambda^2\Phi(x)} \leq \min \{ C_2,\frac 1{\lambda^2\Phi(x)} \}$ by  hypothesis $(\textbf{H}4)$, we find by dominated convergence that $\|\vu_{0,\lambda}\|_{L^2(\Phi\, dx)}=o(\lambda^{\frac{4-d}2})$ and thus $\lim_{ \lambda \rightarrow +\infty} \lambda^2 T_\lambda = + \infty$ .

Thus, if we consider a finite time $T$ and a sequence $\epsilon_k$, we may choose $\lambda$ such that $\lambda^2 T_\lambda >T$ (and such that $\|\vu_{0,\lambda} \|_{L^2(\Phi\, dx)}<\epsilon_0$ if $d=4$);  we have a uniform control of $\vu_{\epsilon,\lambda}$ and of $\vN\otimes\vu_{\epsilon,\lambda}$ on $(0,T_\lambda)$, hence  a uniform control of $\vu_{\epsilon}$ and of $\vN\otimes\vu_{\epsilon}$ on $(0,T)$. We may exhibit  a solution on $(0,T)$ using  the Rellich--Lions theorem by extracting a subsequence $\epsilon_{k_n}$.  A diagonal argument permits then  to obtain a global solution.

Theorem \ref{Thm1} is proved.
\Endproof{}

\section{Proof of Theorem \ref{Thm2} (the case $d=2$).}
 
 In the case of dimension $d=2$, the Navier--Stokes equations are well-posed in $H^1$ and we don't need to mollify the equations. Thus, we may approximate the Navier--Stokes equations with
 
   \begin{equation*}  (NS_\epsilon) 
        \left\{ \begin{matrix} 
            \partial_t \vu_\epsilon  = \Delta \vu_\epsilon  -( \vu_\epsilon \cdot \vN) \vu_\epsilon - \vN p_\epsilon  
            \cr \cr \vN \cdot \vu_\epsilon =0, \phantom{space space} \vu_\epsilon(0,.)=\vu_{0, \epsilon}
        \end{matrix}\right.
    \end{equation*}
    with \begin{equation*}
    \vu_{0, \epsilon} = \mathbb{P}( \phi_\epsilon \vu_0 ) .
\end{equation*}

Then the vorticity  $\omega_\epsilon$ is solution of 
\begin{equation*}   
        \left\{ \begin{matrix} 
            \partial_t \omega_\epsilon  = \Delta \omega_\epsilon  -( \vu_\epsilon \cdot \vN) \omega_\epsilon  
            \cr \cr \vN \cdot \omega_\epsilon =0, \phantom{space space} \omega_\epsilon(0,.)=\omega_{0, \epsilon}
        \end{matrix}\right.
    \end{equation*}
    with \begin{equation*}
    \omega_{0, \epsilon} = \vN\wedge( \phi_\epsilon \vu_0 ) .
\end{equation*}

$\vu_{0,\epsilon}$ belongs to $H^1$, so we know that we have a global  solution  $\vu_\epsilon$. We then just have to prove that, for every finite time $T_0$, we have a uniform control of the norms $\|\omega_\epsilon\|_{L^\infty((0,T_0), L^2(\Phi\, dx))}$ and $\|\vN\omega_\epsilon\|_{L^2 ((0,T_0),L^2(\Phi\, dx))}$.

We can calculate$ \int \partial_t \omega_\epsilon \cdot \omega_\epsilon \Phi \,dx $ so that
\begin{equation*}
    \begin{split}
        \int \frac{ |    \omega_\epsilon  (t,x) | ^2}2  \Phi \, dx &+ \int_0^t  \int   | \vN  \omega_\epsilon   | ^2\, \,    \Phi dx\, ds \\
        &= \int \frac{ |   \omega_{0, \epsilon} (x) | ^2}2    \Phi\, dx  - \int_0^t  \int  \vN(\frac{| \omega_\epsilon |^2} {2} )\, \, \cdot  \vN  \Phi   dx\, ds \\
        & + \int_0^t \int  \frac{| \omega_\epsilon |^2}{2}  \vu_\epsilon  \,  \cdot \vN   \Phi \, dx\, ds.
    \end{split}
\end{equation*}

As
\begin{equation*}
\begin{split}
    \int_0^t \int  \frac{| \omega_\epsilon |^2}{2}  \vu_\epsilon  \,  \cdot \vN   \Phi \, dx\, ds &\leq \int_0^t \| \sqrt{\Phi} \omega_\epsilon \|_{L^{\frac{8}{3}} }^2 \|\sqrt{\Phi} \vu_\epsilon\|_{L^4} \\
    & \leq \int_0^t (\| \sqrt{\Phi} \omega_\epsilon \|_{L^2 }^{3/4} \| \vN (\sqrt{\Phi} \omega_\epsilon) \|_{L^2 }^{1/4} )^2 \|\sqrt{\Phi} \vu_\epsilon\|_{L^4} 
\end{split}
\end{equation*}

we obtain
\begin{equation*}
    \begin{split}
        \| \sqrt{\Phi}  \omega_\epsilon (t) \|_{L^2}^2 +  \int_0^t \| \sqrt{\Phi}  \vN   \omega_\epsilon  \|_{L^2}^2 \leq \| \sqrt{\Phi} \omega_{0, \epsilon} \|_{L^2}^2  + C_\Phi \int_0^t \| \sqrt{\Phi}   \omega_\epsilon  \|_{L^2}^2  (1 + \| \sqrt{\Phi}  \vu_\epsilon  \|_{L^4}^\frac{4}{3}) \, ds 
    \end{split}
\end{equation*}
 
We can conclude that,  for all $T> 0$ and for all $t \in (0,T)$,
\begin{equation*}
    \begin{split}
        \| \sqrt{\Phi}  \omega_{\epsilon} (t) \|_{L^2}^2 +  \int_0^t\| \sqrt{\Phi}  \vN   \omega_{\epsilon}  \|_{L^2}^2
        \leq \| \sqrt{\Phi} \omega_{0, \epsilon} \|_{L^2}^2  e^{ C_\Phi \sup_{\epsilon >0} \int_0^t  (1+ \| \sqrt{\Phi}  \vu_{\epsilon } \|_{L^4})^{\frac{4}{3}} \, ds }
    \end{split}
\end{equation*}
Thus, we have uniform controls on $(0,T)$. \Endproof

\section{Proof of Theorems \ref{Thm3} and \ref{Thm4} (the axisymmetric case)}

\subsection{Axisymmetry.} 

In $\mathbb{R}^3$, we consider the usual coordinates $(x_1,x_2,x_3)$ and  the cylindrical coordinates $(r,\theta,z)$  given by the formulas $x_1 = r \cos \theta, $ $ x_2 = r \sin \theta $ and $x_3=z$.

We denote $(\ve_1 , \ve_2,\ve_3)$ the usual canonical basis
 $$\ve_1=(1, 0, 0),\ \ve_2=(0, 1, 0), \ve_3=(0, 0, 1).$$ We attach to the point $x$ (with $r\neq 0$) another orthonormal basis
\begin{equation*}
    \ve_ r  = \frac{\partial x}{\partial  r } = \cos\theta\,  \ve_1 + \sin\theta\,  \ve_2,   \, \, \,
    \ve_\theta = \frac{1}{ r } \frac{\partial x}{\partial \theta} = -\sin\theta\,  \ve_1 + \cos \theta\,  \ve_2,   \, \, \, 
    \ve_z = \frac{\partial x}{\partial z} =  \ve_3.
\end{equation*}

 For a vector field $\vu  = (u_1, u_2, u_3 ) =u_1 \ve_1 + u_2 \ve_2 + u_3 \ve_3$,  we can see that
$$ \vu = (u_1 \cos \theta + u_2 \sin \theta) \, \ve_ r  + (-u_1 \sin \theta + u_2 \cos \theta) \, \ve_\theta + u_3 \, \ve_z .$$

We will denote $(u_ r , u_\theta, u_z)_p$ the coordinates of $\vu$ in the basis $(\ve_ r , \ve_\theta, \ve_z) $. We will  consider  axially symmetric (axisymmetric) vector fields $\vu$ without swirl and axisymmetric scalar functions $a$, which means that
$$ \vu  = u_r(r,z) \, \ve_r + u_z(r,z) \,  \ve_z \, \, \, \, \,
\text{and} \, \, \, \, \, a = a(r,z). $$

  \subsection{The $H^1$ case.}
  
We will use the following well known results of Ladyzhenskaya \cite{La68, LR16}.     \begin{Proposition}
\label{ple}
Let $\vu_0$ be a divergence free axisymmetric vector field without swirl, such that $\vu_0 $ belongs to $ H^1$.  Then, the following problem
    \begin{equation*}  (NS) 
        \left\{ \begin{matrix} 
            \partial_t \vu= \Delta \vu  -(\vu\cdot \vN)\vu- \vN p 
            \cr \cr \vN\cdot \vu=0, \phantom{space space} \vu(0,.)=\vu_0
        \end{matrix}\right.
    \end{equation*} has a unique solution $\vu\in \mathcal{C}([0,+\infty), H^1)$.
    This solution is axisymmetric 
     without swirl. Moreover,    $\vu, \vN\otimes \vu $ belong to $  L^\infty((0,+\infty),  L^2) $,  and $\vN\otimes  \vu, \Delta \vu $ belong to $L^2((0,+\infty),L^2)$.
     
     If $\vu_0\in H^2$, we have the inequality
       \begin{equation*}
      \int \frac{ | \omega(t) |^2 }{r^2} dx \leq \int \frac{ | \omega_0 |^2 }{r^2}  \leq \| \vN \otimes \omega_0 \|_2^2.
  \end{equation*}
\end{Proposition}

\subsection{A priori controls}
\label{scvst}

Let $\phi \in \mathcal{D}(\mathbb{R}^2)$ be a real-valued radial function which is equal to 1 in a neighborhood of $0$ and let $\phi_\epsilon(x)= \phi(\epsilon(x_1,x_2))$. For $\epsilon\in (0,1]$, let  
\begin{equation*}
    \vu_{0, \epsilon} = \mathbb{P}( \phi_\epsilon \vu_0 ) .
\end{equation*}
Thus, $\vu_{0, \epsilon}$ is a divergence free axisymmetric without swirl vector field which belongs to $H^1$. As we have 
\begin{equation*}
    \omega_{0, \epsilon} = \vN \wedge \vu_{0, \epsilon} = \vN \wedge (\phi_\epsilon \vu_0) = \phi_\epsilon \omega_0 + \epsilon (\vN \phi)(\epsilon x) \wedge \vu_0,
\end{equation*}
using $\Phi\in\mathcal{A}_2$ and $\vert \epsilon \vN\phi(\epsilon x)\vert \leq C \frac 1 r \mathds{1}_{  r\geq \frac 1 {C\epsilon}}\leq C'  \mathds{1}_{  r \geq \frac 1 {C\epsilon}}\sqrt\Phi$, we can see that
\begin{equation*}
    \lim_{\varepsilon \to 0} \| \vu_0 - \vu_{0, \epsilon} \|_{L^2 (\Phi\,  dx)} + \| \omega_0 - \omega_{0, \epsilon} \|_{L^2 (\Psi \, dx)} =0.
\end{equation*}
 
Let $\vu_\epsilon$ be the global solution of the problem 
    \begin{equation*}  (NS_\epsilon) 
        \left\{ \begin{matrix} 
            \partial_t \vu_\epsilon  = \Delta \vu_\epsilon  -(\vu_\epsilon \cdot \vN) \vu_\epsilon - \vN p_\epsilon
            \cr \cr \vN\cdot \vu_\epsilon=0, \phantom{space space} \vu_\epsilon(0,.)=\vu_{0, \epsilon}
        \end{matrix}\right.
    \end{equation*}
given by the Proposition \ref{ple}. We denote $ \omega_\epsilon = \vN\wedge\,  \vu_\epsilon$, then
\begin{equation}
\label{vue}
    \partial_t \vu_\epsilon = \Delta \vu_\epsilon + ( \vu_\epsilon \cdot \vN)\vu_\epsilon - \vN p_\epsilon 
\end{equation} 
and
\begin{equation}
\label{we}
    \partial_t \omega_\epsilon = \Delta \omega_\epsilon + ( \omega_\epsilon \cdot \vN) \vu_\epsilon - (\vu_\epsilon \cdot \vN) \omega_\epsilon 
\end{equation} 

As $\sqrt{\Psi} \omega_\epsilon \in L^2 H^1$ (because $\sqrt{\Psi},\vN \sqrt{\Psi} \in L^\infty$) and $\sqrt{\Psi} \partial_t \omega_\epsilon \in L^2 H^{-1}$, we can calculate$ \int \partial_t \omega_\epsilon \cdot \omega_\epsilon \Psi \,dx $ using \eqref{we} so that
\begin{equation*}
    \begin{split}
        \int \frac{ |    \omega_\epsilon  (t,x) | ^2}2  \Psi \, dx &+ \int_0^t  \int   | \vN  \otimes\omega_\epsilon   | ^2\, \,    \Psi dx\, ds \\
        =& \int \frac{ |   \omega_{0, \epsilon} (x) | ^2}2    \Psi\, dx  - \int_0^t  \int  \vN(\frac{| \omega_\epsilon |^2} {2} )\, \, \cdot  \vN  \Psi   dx\, ds \\
        & + \int_0^t \int  \frac{| \omega_\epsilon |^2}{2}  \vu_\epsilon  \,  \cdot \vN   \Psi \,  - (\omega_\epsilon\cdot\vu_\epsilon)\omega_\epsilon\cdot\vN\Psi\, dx\\
        & - \int_0^t \int   (( \omega_\epsilon  \cdot \vN)  \omega_\epsilon  ) \cdot  \vu_\epsilon  \,   \,  \Psi  \, dx\, ds   
        \\    \leq& \int \frac{ |   \omega_{0, \epsilon} (x) | ^2}2    \Psi\, dx  +\frac 1 8 \int_0^t  \int   | \vN  \otimes\omega_\epsilon   | ^2\, \,    \Psi dx\, ds + C \int_0^t \|\sqrt\Psi\, \omega_\epsilon\|_2^2\, ds    \\
        & + C  \int_0^t  \|\sqrt\Psi \,  \omega_\epsilon\|_2   \|\sqrt\Psi \,  \omega_\epsilon\|_6     \| \sqrt\Phi \vu_\epsilon\|_3\, ds \\
        & - \int_0^t \int   (( \omega_\epsilon  \cdot \vN)  \omega_\epsilon  ) \cdot  \vu_\epsilon  \,   \,  \Psi  \, dx\, ds \\    \leq& \int \frac{ |   \omega_{0, \epsilon} (x) | ^2}2    \Psi\, dx  +\frac 1 4 \int_0^t  \int   | \vN  \otimes\omega_\epsilon   | ^2\, \,    \Psi dx\, ds + C \int_0^t \|\sqrt\Psi\, \omega_\epsilon\|_2^2\, ds    \\
        & + C'   \int_0^t  \|\sqrt\Psi \,  \omega_\epsilon\|_2^2        (\| \sqrt\Phi \vu_\epsilon\|_3+  (\| \sqrt\Phi \vu_\epsilon\|_3^{4/3}) \, ds \\
        & - \int_0^t \int   (( \omega_\epsilon  \cdot \vN)  \omega_\epsilon  ) \cdot  \vu_\epsilon  \,   \,  \Psi  \, dx\, ds   \end{split}
\end{equation*}

As $ \omega_\epsilon  = \omega_{\epsilon, \theta}  \, \ve_\theta $, we have $$  \omega_\epsilon  \cdot \vN  \omega_\epsilon  = -\frac{ \omega_{\epsilon, \theta} ^2 }{r}  \, \ve_r  .$$ 

In order to control $ \vu_\epsilon  \cdot (  \omega_\epsilon  \cdot  \vN  \omega_\epsilon )$, we split the domain of  integration in a domain where $r$ is small and a domain where $r$ is large. The support of $\phi_1$   is contained in $\{x\ /\ r<R\}$ for some $R>0\}$, and the support of $1-\phi_1$ is contained in $\{x\ /\ r>R_0\}$ for some $R_0>0\}$. We have
$$ \inf_{r<R} \Phi(x)= \inf_{\sqrt {x_1^2+x_2^2}<R}\Phi(x_1,x_2,0)>0$$ and similarly
$$ \inf_{r<R} \Psi(x)= \inf_{\sqrt {x_1^2+x_2^2}<R}\Psi(x_1,x_2,0)>0.$$  
On the other hand, we have
$$ \inf_{r>R_0} r^2 \Phi(x)=\inf_{\sqrt {x_1^2+x_2^2}>R_0} (x_1^2+x_2^2) \Phi(x_1,x_2,0)\geq \inf_{\vert x\vert>R_0} \vert x\vert^2 \Phi(x)>0.$$
We then write:
\begin{equation*}
    \begin{split}
     - \int_0^t \int &   (( \omega_\epsilon  \cdot \vN)  \omega_\epsilon  ) \cdot  \vu_\epsilon   \,   \,  \Psi  \, dx\, ds \\ 
        = &  \int_0^t \int  \phi_1( \, ( \omega_\epsilon  \cdot \vN)  \vu_\epsilon  ) \cdot  \omega_\epsilon  \, ) \,  \Psi  \, dx\, ds + \int_0^t \int  ( \omega_\epsilon  \cdot  \vu_\epsilon ) ( \omega_\epsilon  \cdot \vN \phi_1)   \Psi  \, dx\, ds   \\
        & +\int_0^t  \int  \phi_1 (  \omega_\epsilon  \cdot  \vu_\epsilon  )   \omega_\epsilon  \cdot  \vN   \Psi   dx \, ds \\
        & - \int_0^t  \int (1-\phi_1) (  \vu_\epsilon  \cdot (  \omega_\epsilon  \cdot  \vN  \omega_\epsilon ) ) \Psi dx \, ds  \\ 
       \leq &C   \int_0^t \int    \vert \omega_\epsilon\vert^2  \vert  \vN \otimes  \vu_\epsilon \vert \,  \Psi^{3/2}  \, dx\, ds + C \int_0^t \int    \vert \omega_\epsilon\vert^2  \vert    \vu_\epsilon \vert \,  \sqrt\Phi\,    \Psi  \, dx\, ds .   \end{split}
\end{equation*}
 As $\Psi\in\mathcal{A}_2$, we have $\|\sqrt\Psi \vN\otimes\vu_\epsilon\|_2\approx \|\sqrt\Psi\omega_\epsilon\|_2$; moreover, $$\|\vN\otimes(\sqrt\Phi\vu_\epsilon)\|_2\leq C(\|\sqrt\Phi \vu_\epsilon\|_2+\|\sqrt\Psi \omega_\epsilon\|_2)$$ and $$ \|\vN\otimes(\sqrt\Psi\omega_\epsilon)\|_2\leq C(\|\sqrt\Psi \omega_\epsilon\|_2+\|\sqrt\Psi  \vN\otimes\omega_\epsilon\|_2),$$ and thus we get
\begin{equation*}
    \begin{split}
     - \int_0^t \int &    (( \omega_\epsilon  \cdot \vN)  \omega_\epsilon  ) \cdot  \vu_\epsilon   \,   \,  \Psi  \, dx\, ds \\ 
        \leq &  C \int_0^t \| \sqrt{\Psi} \vN \otimes
          \vu_\epsilon  \|_{L^2} \| \sqrt{\Psi}  \omega_\epsilon  \|_{L^3} \| \sqrt{\Psi}  \omega_\epsilon  \|_{L^6} \, ds \\
          & + C \int_0^t \| \sqrt{\Phi}   \vu_\epsilon  \|_{L^6} \| \sqrt{\Psi}  \omega_\epsilon  \|_{L^3} \| \sqrt{\Psi}  \omega_\epsilon  \|_{L^2} \, ds\\
        \leq & C' \int_0^t \| \sqrt{\Psi}    \omega_\epsilon  \|_{L^2}^{\frac{3}{2}} (\| \sqrt{\Psi}  \omega_\epsilon  \|_{L^2}  + \| \sqrt{\Psi}  \vN\otimes \omega_\epsilon  \|_{L^2} )^{\frac{3}{2}} \, ds\\
        & +  C'  \int_0^t \| \sqrt{\Phi}  \vu_\epsilon  \|_{L^2}  \| \sqrt{\Psi}    \omega_\epsilon  \|_{L^2}^{\frac{3}{2}} (\| \sqrt{\Psi}  \omega_\epsilon  \|_{L^2}  + \| \sqrt{\Psi}  \vN \otimes\omega_\epsilon  \|_{L^2} )^{\frac{1}{2}} \, ds\\
       \leq &   C'' \int_0^t ( \| \sqrt{\Phi}  \vu_\epsilon  \|_{2} +\| \sqrt{\Phi}  \vu_\epsilon  \|_{2} ^{4/3}) \|\sqrt\Psi\omega_\epsilon\|_2^2+ \| \sqrt{\Psi}  \omega_\epsilon  \|_{2}^{3} + \| \sqrt{\Psi}  \omega_\epsilon  \|_{2}^6 \, ds\\&+ \frac{1}{8} \int_0^t \| \sqrt{\Psi}  \vN \otimes  \omega_\epsilon  \|_{2}^2\, ds
    \end{split}
\end{equation*}

 We finally find that
\begin{equation}
\label{coe1}
    \begin{split}
        &\| \sqrt{\Psi}  \omega_\epsilon (t) \|_{L^2}^2 + \int_0^t \| \sqrt{\Psi}  \vN   \otimes\omega_\epsilon  \|_{L^2}^2\, ds \\
        & \leq \| \sqrt{\Psi} \omega_{0, \epsilon} \|_{L^2}^2  + C \int (1+\| \sqrt\Phi \vu_\epsilon\|_3+  (\| \sqrt\Phi \vu_\epsilon\|_3^{4/3})  \|\sqrt\Psi \omega_\epsilon\|_2^2   \, ds \\& +  C \int_0^t ( \| \sqrt{\Phi}  \vu_\epsilon  \|_{2}\ +\| \sqrt{\Phi}  \vu_\epsilon  \|_{2} ^{4/3}) \|\sqrt\Psi\omega_\epsilon\|_2^2+ \| \sqrt{\Psi}  \omega_\epsilon  \|_{2}^{3} + \| \sqrt{\Psi}  \omega_\epsilon  \|_{2}^6 \, ds \\
          & \leq \| \sqrt{\Psi} \omega_{0, \epsilon} \|_{L^2}^2  \\& +  C' \int_0^t (1+ \| \sqrt{\Phi}  \vu_\epsilon  \|_{2}\ +\| \sqrt{\Phi}  \vu_\epsilon  \|_{2} ^{4/3}) \|\sqrt\Psi\omega_\epsilon\|_2^2+  \| \sqrt{\Psi}  \omega_\epsilon  \|_{2}^6 \, ds 
    \end{split}
\end{equation}

We already know that $\| \sqrt{\Phi}  \vu_\epsilon (t) \|_{L^2}$ remains bounded (independently of $\epsilon$) on every bounded interval, so that we may again use Lemma \ref{lem5} and control 
$  \sup_{0\leq t\leq T_0} \|\ \omega_\epsilon (t,.)\|_{L^2(\Psi dx)}^2 +  \int_0^{T_0} \|\vN \omega_\epsilon \|_{L^2(\Psi dx)}^2\, ds$  for some $T_0$, where both $T_0$ and the control don't depend on $\epsilon$.  The control is then transferred to the limit $\omega$ since $\omega=\lim \omega_{\epsilon_k}=\lim \vN\wedge \vu_{\epsilon_k}$.
This proves local existence of a regular solution and Theorem \ref{Thm3} is proved.



\subsection{The case of a very regular initial value.}

We present a result apparently more restrictive that our main Theorem (Theorem \ref{Thm4}), but we will see that it implies almost directly our main Theorem.

\begin{Proposition}
\label{lge}

Let $ \Phi $ be a weight satisfying $(\textbf{H}1)-(\textbf{H}4) $.  Assume moreover that $\Phi $ depends only on $r=\sqrt{x_1^2+x_2^2}$.  Let $\Psi$ be another continuous weight (that  depends only on $r$)   such that  $\Phi \leq \Psi\leq 1$, $\Psi\in\mathcal{A}_2$ and       there exists $C_1>0$ such that $$ \vert \vN \Psi  \vert \leq C_1  \sqrt\Phi  \Psi \text{ and }  \vert\Delta \Psi\vert\leq C_1 \Phi\Psi.$$

Let $\vu_0$ be a divergence free axisymmetric vector field without swirl, such that $\vu_0,$ belongs to $ L^2 (\Phi dx) $,  $\vN\otimes \vu_0 $ and $\Delta  \vu_0$ belong to $ L^2 (\Psi dx) $.   Then there exists a global solution $\vu$ of the problem
    \begin{equation*}  (NS) 
        \left\{ \begin{matrix} 
            \partial_t \vu= \Delta \vu  -(\vu\cdot \vN)\vu- \vN p 
            \cr \cr \vN\cdot \vu=0, \phantom{space space} \vu(0,.)=\vu_0
        \end{matrix}\right.
    \end{equation*}
such that
 \begin{itemize}
    \item $\vu $ is axisymmetric without swirl, $\vu$    belongs to $ L^\infty ((0,T), L^2 (\Phi\,  dx) ) $,  $\vN\otimes \vu $ belong to $ L^\infty ((0,T), L^2 (\Psi\,  dx) ) $ and $ \Delta \vu $ belongs to $L^2((0,T),L^2 (\Psi\,  dx) )$, for all $T>0$,
     
    \item the maps $t\in [0,+\infty)\mapsto \vu(t,.)$ and $t\in [0,+\infty)\mapsto \vN \otimes\vu(t,.)$ are weakly continuous from $[0,+\infty)$ to $L^2 (\Phi\,  dx) $ and to $L^2(\Psi\, dx)$ respectively, and are strongly continuous at $t=0$.
         
 \end{itemize}
 
\end{Proposition}

\pv

  Ladyzhenskaya's inequality for axisymmetric fields with no swirl    (Proposition \ref{ple}) gives
\begin{equation}
    \label{vr2}
    \begin{split}
        \int \frac{| \omega_\epsilon (t)|^2}{r^2} dx \leq \int \frac{|\omega_{0, \epsilon}|^2}{r^2} dx.
    \end{split}
\end{equation}

As we have 
\begin{equation*}
    \partial_i \omega_{0, \epsilon}   = \phi_\epsilon \partial_i \omega_0 + \epsilon \partial_i \phi(\epsilon x)  \omega_0  + \epsilon (\vN \phi)(\epsilon x) \wedge \partial_i \vu_0 + \epsilon^2  (\vN \partial_i \phi)(\epsilon x) \wedge \vu_0,
\end{equation*}
we can see that
\begin{equation*}
    \lim_{\epsilon \to 0 } \| \vN \otimes\omega_{0, \epsilon} - \vN\otimes \omega_{0} \|_{L^2 (\Psi\, dx)} = 0.
\end{equation*}
As
\begin{equation*}
    \int \frac{|\omega_{0, \epsilon} - \omega_{0} |^2}{r^2} dx \leq C (\int_{ 0 < r < 1 } |\vN\otimes  \omega_{0, \epsilon} - \vN \otimes\omega_{0} |^2 \Psi\,  dx + \int_{ 1 < r < + \infty } |\omega_{0, \epsilon} - \omega_{0} |^2 \Psi\, dx),
\end{equation*}
we also have
\begin{equation*}
    \lim_{\epsilon \to 0 } \int \frac{|\omega_{0, \epsilon} - \omega_{0} |^2}{r^2} dx =0.
\end{equation*}

We know that
\begin{equation*}
    \begin{split}
        &\int \frac{ |    \omega_\epsilon  (t,x) | ^2}2  \Psi \, dx + \int_0^t  \int   | \vN  \otimes\omega_\epsilon   | ^2\, \,    \Psi dx\, ds \\
        &= \int \frac{ |   \omega_{0, \epsilon} (x) | ^2}2    \Psi\, dx  - \int_0^t  \int  \vN(\frac{| \omega_\epsilon |^2} {2} )\, \, \cdot  \vN  \Psi   dx\, ds \\
        & + \int_0^t \int  \frac{| \omega_\epsilon |^2}{2}  \vu_\epsilon  \,  \cdot \vN   \Psi \, dx\, ds \\
        & - \int_0^t \int  ( \omega_\epsilon  \cdot   \vu_\epsilon  )   \omega_\epsilon  \, \cdot  \vN \Psi  \, dx\, ds - \int_0^t \int  \vu_\epsilon (  \omega_\epsilon \cdot \vN \omega_\epsilon  ) \,  \Psi  \, dx\, ds  
    \end{split}
\end{equation*}
which implies
\begin{equation*}
    \begin{split}
        &\| \sqrt{\Psi}  \omega_\epsilon (t) \|_{L^2}^2 + 2 \int_0^t \| \sqrt{\Psi}  \vN   \omega_\epsilon  \|_{L^2}^2 \\
        & \leq \| \sqrt{\Psi} \omega_{0, \epsilon} \|_{L^2}^2 + 2 \int_0^t \| \sqrt{\Psi}  \omega_\epsilon  \|_{L^2} \| \sqrt{\Psi} \vN  \omega_\epsilon  \|_{L^2} \\
        & + \int_0^t \| \sqrt{\Phi}   \vu_\epsilon  \|_{L^3} \| \sqrt{\Psi}  \omega_\epsilon  \|_{L^3}^2 \\
        & + \int_0^t  \frac{1}{r} |\vu_{r, \epsilon}| |\omega_\epsilon|^2 \Psi\, dx \, ds 
    \end{split}
\end{equation*}
Furthermore, we have
 \begin{equation*}
     \int_0^t  \int \frac{1-\phi_1(x)}{r} |\vu_{r, \epsilon}| |\omega_\epsilon|^2 \Psi\,  dx \, ds \leq \int_0^t \| \sqrt{\Phi}   \vu_\epsilon  \|_{L^3} \| \sqrt{\Psi}  \omega_\epsilon  \|_{L^3}^2
 \end{equation*}
 and
  \begin{equation*}
     \int_0^t  \int \frac{ \phi_1(x)}{r} |\vu_{\epsilon,r}| |\omega_\epsilon|^2 dx \, ds \leq C \int_0^t \| \frac{\omega _\epsilon}{r}  \|_{L^2} \| \sqrt{\Psi}  \vu_\epsilon  \|_{L^\infty} \| \sqrt{\Psi}  \omega_\epsilon  \|_{L^2},
 \end{equation*}
where
\begin{equation*}
    \begin{split}
        \| \frac{\omega _\epsilon}{r}  \|_{L^2} \leq C \| \frac{\omega _{0, \epsilon} }{r}  \|_{L^2} &\leq C ( \| \sqrt{\Psi} \omega _{0, \epsilon}   \|_{L^2}+ \| \sqrt{\Psi} \vN\otimes \omega _{0, \epsilon}    \|_{L^2} ) \\
        & \leq C' ( \| \sqrt{\Phi} \vu _{0}   \|_{L^2} + \| \sqrt{\Psi} \omega _{0}   \|_{L^2}+\| \sqrt{\Psi} \vN \otimes\omega _{0}    \|_{L^2} )
    \end{split}
\end{equation*} 
 and
\begin{equation*}
    \begin{split}
    \| \sqrt{\Psi}  \vu_\epsilon  \|_{L^\infty}^2 &\leq C \| \vN\otimes (\sqrt\Psi \vu_\epsilon) \|_{2} \| \Delta (\sqrt\Psi \vu_\epsilon) \|_{2} \\
    &\leq C' (\| \sqrt{\Phi}  \vu_\epsilon  \|_{L^2} + \| \sqrt{\Psi}  \omega_\epsilon  \|_{L^2} + \| \sqrt{\Psi} \vN\otimes \omega_\epsilon  \|_{L^2})^2.
    \end{split}
\end{equation*}

Then, if we denote $A_0 = \| \sqrt{\Phi}  \vu_{0}  \|_{L^2} + \| \sqrt{\Psi}  \omega_0  \|_{L^2} + \| \sqrt{\Psi} \vN\otimes \omega_0 \|_{L^2}$, we get
\begin{equation*}
    \begin{split}
        \| \sqrt{\Psi}  \omega_\epsilon (t) \|_{L^2}^2 &+  \int_0^t \| \sqrt{\Psi}  \vN  \otimes \omega_\epsilon  \|_{L^2}^2 \\
        \leq& \| \sqrt{\Psi} \omega_{0, \epsilon} \|_{L^2}^2 + C \int_0^t \| \sqrt{\Phi}  \vu_\epsilon  \|_{L^2}^2 \\
        & + C_\Phi \int_0^t \| \sqrt{\Psi}   \omega_\epsilon  \|_{L^2}^2 (1 + A_0 + A_0^2 + \| \sqrt{\Phi}  \vu_\epsilon  \|_{L^3} + \| \sqrt{\Phi}  \vu_\epsilon  \|_{L^3}^2) \, ds 
    \end{split}
\end{equation*}

So, we conclude that, for all $T> 0$ and for all $t \in (0,T)$,
\begin{equation*}
    \begin{split}
        \| & \sqrt{\Psi}  \omega_{\epsilon} (t) \|_{L^2}^2 +  \int_0^t\| \sqrt{\Psi}  \vN   \otimes\omega_{\epsilon}  \|_{L^2}^2 \\
        \leq& (\| \sqrt{\Psi} \omega_{0, \epsilon} \|_{L^2}^2 + C_\Phi \sup_{\epsilon >0}   \int_0^T \| \sqrt{\Phi}  \vu_{\epsilon}  \|_{L^2}^2) e^{ C_\Phi \sup_{\epsilon >0} \int_0^t (1  + A_0^2 + \| \sqrt{\Phi}  \vu_{\epsilon}  \|_{L^3} + \| \sqrt{\Phi}  \vu_{\epsilon } \|_{L^3}^2) \, ds }
    \end{split}
\end{equation*}
Thus, we can obtain a solution on $(0,T)$ using the Aubin--Lions Theorem and finish with a  diagonal argument to get a global  solution.
\Endproof{}



\subsection{End of the proof.}

We begin by consider a local solution $\vv$ on $(0,T_0)$ with initial value $\vu_0$ given by Theorem \ref{Thm3}, which is  continuous from $(0, T_0)$ to $\mathcal{D}'$. We take $T_1 \in (0,T_0)$ such that $\vN\otimes (\vN \wedge \vv)(T_1,.) \in L^2(\Phi dx)$. We consider a solution $\textbf{w}$ on $(T_1, + \infty)$ and initial value $\vv(T_1)$ given by Proposition  \ref{lge}. Our global solution is defined as $\vu=\vv $ on $(0,T_1)$ and  $\vu=\vw$ on $(T_1, + \infty)$. \Endproof{}


\end{document}